\input amstex
\documentstyle{amsppt}
\NoBlackBoxes
\TagsOnRight
\NoRunningHeads
\hcorrection{2cm}

\loadbold
\loadeusm
\loadeufm

\def\qed{\hphantom{a}\hfill $\square$}
\def\ee{\hphantom{a}\hfill $\lozenge$}
\def\Ra{\Rightarrow}

\def\LRa{\Leftrightarrow}
\def\alb{\allowmathbreak}
\def\mb{\mathbreak}
\def\nmb{\nomathbreak}

\def\supp{\operatorname{supp}}

\def\Aut{\operatorname{Aut}}

\def\<{\langle}
\def\>{\rangle}

\def\slim{\operatorname{s-lim}}
\def\wlim{\operatorname{w-lim}}
\def\w*lim{\operatorname{w*-lim}}

\def\vep{\varepsilon}
\def\vp{\varphi}
\def\vt{\vartheta}

\def\sub{\subseteq}

\def\ov{\overline}
\def\r{\negthickspace\restriction\negthickspace}

\def\G{\varGamma}
\def\g{\gamma}
\def\om{\omega}
\def\Om{\varOmega}
\def\be{\beta}
\def\al{\alpha}
\def\Lam{\Lambda}
\def\lam{\lambda}
\def\del{\delta}
\def\Del{\Delta}
\def\iy{\infty}
\def\si{\sigma}
\def\k{\kappa}
\def\ti{\times}
\def\oti{\otimes}
\def\inl{\int\limits}
\def\di{\diamond}
\def\tr{\operatorname{tr}}
\topmatter
\title{The Choquet-Deny equation in a Banach space}
\endtitle
\author{\smc Wojciech Jaworski* and Matthias Neufang*}
\vskip1pt {\rm School of Mathematics and 
Statistics \\ Carleton University\\
Ottawa, Ontario, Canada K1S 5B6}
\vskip2pt
{\rm e-mails:\ {\it wjaworsk\@math.carleton.ca\/} and {\it 
mneufang\@math.carleton.ca}}
\endauthor 
\leftheadtext{}
\rightheadtext{}
\abstract Let $G$ be a locally compact group and $\pi$ a representation of 
$G$ by weakly* continuous isometries acting in a dual Banach space $E$. 
Given a 
probability measure $\mu$ on $G$ we study the Choquet-Deny equation 
$\pi(\mu)x=x$, $x\in E$. We prove that the solutions of this equation  
form the range of a projection of norm 1 and can be represented by means of a 
``Poisson formula'' on the same boundary space that is used to represent the 
bounded harmonic functions of the random walk of law $\mu$. The relation 
between the space of solutions of the Choquet-Deny equation in $E$ and the 
space of bounded harmonic functions can be understood in terms of a 
construction resembling the $W^*$-crossed product and coinciding precisely 
with the crossed product in the special case of the Choquet-Deny equation in 
the space $E=B(L^2(G))$ of bounded linear operators on $L^2(G)$. Other 
general properties of the Choquet-Deny equation in a Banach space are also 
discussed. 
\medpagebreak

\endabstract
\endtopmatter
\document
\footnote""{*\ Supported by NSERC Grants.}
\head{\bf 1. Introduction} \endhead
The classical Choquet-Deny theorem asserts that when $\mu$ is a regular 
probability measure on a locally compact abelian group $G$, then 
every bounded continuous function $h\:G\to\Bbb C$ which satisfies 
$$
h(g)=\int_G h(gg')\,\mu(dg')\tag{1.1}
$$
for every $g\in G$, is necessarily constant on the cosets of the smallest 
closed subgroup, $G_\mu$, containing the support of $\mu$. The theorem can be 
readily seen to be equivalent to the statement that whenever 
$\pi$ is a representation of $G$ by weakly continuous isomorphisms of a 
locally convex space $E$, where for 
every $x\in E$ and $x^*\in E^*$ the function $G\ni g\to \<\pi(g)x,x^*\>\in 
\Bbb C$ is bounded and continuous, then every vector $x\in E$ which satisfies 
$$
x= \int_G \pi(g)x\medspace\mu(dg),\tag{1.2}
$$
is necessarily fixed by every $g\in G_\mu$, i.e., $\pi(g)x=x$ for every 
$g\in G_\mu$. 

The Choquet-Deny theorem is known to remain true for some nonabelian 
groups too; however, it is not true for all groups. When the theorem 
fails then Eqs (1.1) and (1.2) will have other solutions in 
addition to those described above. The goal of this article is to understand 
the connection between the solutions of the classical 
{\it Choquet-Deny equation\/} (1.1) and its functional analytic 
counterpart (1.2), and to uncover some general properties of the spaces of 
such solutions. We will be concerned mainly with representations that are 
adjoints of strongly continuous representations by isometries in separable 
Banach spaces and our group $G$ will, for the most part, be locally compact 
and second countable. 

In general, a bounded Borel function $h\:G\to \Bbb C$ satisfying (1.1) is 
called a {\it $\mu$-harmonic} function. $\mu$-harmonic functions have their 
origin in probability theory. They are a special case of the harmonic 
functions of a Markov chain and, as such, have been extensively studied. 
When $G$ is locally compact one usually considers the bounded 
$\mu$-harmonic functions as elements of $L^\iy(G)$. They are then solutions 
of Eq.\,(1.2) where $\pi$ is now the right regular representation in 
$L^\iy(G)$. In fact, 
the right hand sides of Eqs (1.1) and (1.2) define a positive weakly* 
continuous contraction $\pi(\mu)\:L^\iy(G)\to L^\iy(G)$ which commutes with 
the left regular representation $\pi_l$ of $G$ in $L^\iy(G)$ 
and can be viewed as the average $\int_G \pi(g)\,\mu(dg)$ of the 
right regular representation. It is then immediate that the space 
$\Cal H_\mu\sub L^\iy(G)$ of bounded $\mu$-harmonic functions is a 
weakly* closed self-adjoint 
subspace of $L^\iy(G)$ containing constants and invariant under $\pi_l$. It 
is a much deeper result, involving the theory of martingales, that bounded 
$\mu$-harmonic functions can be represented, by means of a ``Poisson 
formula'', as bounded Borel functions on a certain ``boundary space''. As a 
result, $\Cal H_\mu$ turns out to be isomorphic to an $L^\iy$-space and so 
has canonically the structure of an abelian $W^*$-algebra. More precisely, 
for each $h_1,h_2\in \Cal H_\mu$ the sequence $\pi(\mu^n)(h_1h_2)$ 
converges almost everywhere (hence, also weakly*), and the formula 
$$
h_1\di h_2=\w*lim_{n\to\iy}\pi(\mu^n)(h_1h_2) \tag{1.3}
$$ 
defines a product $\di$ on $\Cal H_\mu$ under which $\Cal H_\mu$ is an 
abelian $W^*$-algebra. The product $\di$ coincides with the ordinary product 
in $L^\iy(G)$ if and only if the Choquet-Deny theorem holds for the 
particular measure $\mu$. It is also remarkable that $\Cal H_\mu$ is the 
range of a projection $K\: L^\iy(G)\to\Cal H_\mu$ of norm 1, 
which commutes with the left regular representation. 

We will show that the solutions of the Choquet-Deny equation (1.2) in 
a dual Banach space also form the range of a projection of norm 
1 and can be represented by means of a Poisson formula on the same boundary 
space that is used to represent the classical bounded $\mu$-harmonic 
functions. It will follow that the relation 
between the space of such solutions and the space $\Cal H_\mu$ 
of the classical $\mu$-harmonic functions can be understood in terms of a 
construction closely resembling the well known construction of the $W^*$-
crossed product, and coinciding precisely with the latter in the 
special case of the Choquet-Deny equation in the space $B(L^2(G))$ of bounded 
linear operators on $L^2(G)$. We will also see that the space of solutions of 
the Choquet-Deny equation in a $W^*$-algebra has itself canonically the 
structure of a $W^*$-algebra. 

\head{\bf 2. Some special cases}\endhead

Let $E$ be a Banach space which is the dual of a Banach space $E_*$. The 
duality between $E_*$ and $E$ will be written $\<x_*\,,x\>$, $x_*\in E_*$, 
$x\in E$. Let $\pi_*$ be a strongly continuous representation of a 
locally compact group $G$ by isometries in $E_*$ and let $\pi$ be the adjoint 
of $\pi_*$, i.e., the representation in $E$ by weakly* continuous isometries 
given by $\pi(g)=[\pi_*(g^{-1})]^*$. As is well known, the representations 
$\pi_*$ and $\pi$ can be extended to representations of the measure algebra 
$M(G)$ in $E_*$ and $E$, resp.. Given $\si\in M(G)$, $\pi_*(\si)$ and 
$\pi(\si)$ are the bounded linear operators which satisfy 
$$\gathered 
\<\pi_*(\si)x_*\,,x\>= \int_G \<\pi_*(g)x_*\,,x\>\,\si(dg)\\
\text{and} \hphantom{ccccccccccccccccccccccccc\<\pi_*(\si)x_*\,,x\>= \int_G 
\<\pi_*(g)x_*\,,x\>\,\si(dg)} \\
\<x_*\,,\pi(\si)x\>= \int_G \<x_*\,,\pi(g)x\>\,\si(dg)\endgathered\tag{2.1}
$$
for all $x_*\in E_*$ and $x\in E$. Clearly, $\|\pi_*(\si)\|, \|\pi(\si)\|
\le\|\si\|$ and 
$\pi(\si)=[\pi_*(\tilde\si)]^*$ where $\tilde\si$ denotes the 
measure $\tilde\si(A)=\si(A^{-1})$. The Choquet-Deny equations in $E_*$ and 
$E$ become 
$$
\pi_*(\mu)x_*=x_*\quad\text{and}\quad\pi(\mu)x=x\tag{2.2}
$$
where $\mu\in M_1(G)$, the subset of probability measures in $M(G)$. 
Solutions of Eqs\ (2.2) will be referred to as {\it $\mu$-harmonic 
vectors\/} and we will denote by $\Cal H_{\mu\pi_*}\sub E_*$ and $\Cal 
H_{\mu\pi}\sub E$ the spaces of such vectors. We note that 
$\Cal H_{\mu\pi_*}$ and $\Cal H_{\mu\pi}$ always include all those vectors 
that are fixed by every element of the subgroup $G_\mu$. These are the 
{\it trivial\/} solutions of the Choquet-Deny equation. 

In the classical case we have $E=L^\iy(G)$, $E_*=L^1(G)$, and $\pi$ and 
$\pi_*$ are the usual right regular representations in $L^\iy(G)$ and 
$L^1(G)$. While the Choquet-Deny theorem is not true in general, it belongs 
to the folklore of the subject that the Choquet-Deny equation in $L^1(G)$, 
and, more generally, in $L^p(G)$, $1\le p<\iy$, as well as the Choquet-Deny 
equations in $C_0(G)$ and $M(G)$ admit only trivial solutions, no matter what 
$G$ and $\mu$ are. This fact, as well as a more general conclusion in the 
same direction, can be easily deduced from the following fundamental result 
due to Mukherjea [30] and Derriennic [9]. 

\proclaim{Lemma 2.1} If $G_\mu$ is not compact then the convolution 
powers $\mu^n$ converge to zero in the weak* topology of $M(G)$. 
\endproclaim 

\proclaim{Proposition 2.2} Let $\pi$ be a representation of $G$ by weakly 
continuous isomorphisms of a locally 
convex space $E$ and suppose that $\pi$ vanishes at infinity, i.e., for every 
$x\in E$ and $x^*\in E^*$ the function $G\ni g\to \<\pi(g)x\,,x^*\>$ 
belongs to $C_0(G)$. Let $v$ be a solution of the Choquet-Deny equation in 
$E$, i.e., 
$$
\<v\,,x^*\>=\int_G\<\pi(g)v\,,x^*\>\,\mu(dg)\tag{2.3}
$$
for every $x^*\in E^*$. Then $\pi(g)v=v$ for every $g\in G_\mu$\,\rom; if 
$G_\mu$ is not compact then $v=0$.
\endproclaim 
\demo{Proof}If $G_\mu$ is compact, consider the restriction of $\pi$ to 
$G_\mu$ and consider $\mu$ as a measure on $G_\mu$. Since the Choquet-Deny 
theorem holds for compact groups, the result follows without difficulty 
because the function $g\to\<\pi(g)v\,,x^*\>$ is $\mu$-harmonic. If 
$G_\mu$ is not compact, observe that Eq.\,(2.3) holds with 
$\mu$ replaced by $\mu^n$; then use Lemma 2.1. 
\qed\enddemo

\proclaim{Corollary 2.3} Let $\pi_*$ and $\pi$ denote the right regular 
representations of $G$ in $C_0(G)$ and $M(G)$, resp., and let $\mu\in 
M_1(G)$. Then $h\in\Cal H_{\mu\pi_*}$ if and only if $h$ is constant on the 
left cosets of $G_\mu$\,\rom; $\si\in\Cal H_{\mu\pi}$ if and only if $\si 
*\del_g=\si$ for every $g\in G_\mu$. If $G_\mu$ is not compact then $\Cal 
H_{\mu\pi_*} =\{0\}$ and $\Cal H_{\mu\pi}=\{0\}$. 
\endproclaim 

\proclaim{Corollary 2.4} Let $\pi_*$ be the right regular representation of 
$G$ in $L^p(G)$ where $1\le p <\iy$. Then $f\in \Cal H_{\mu\pi_*}$ if and 
only if for every $g\in G_\mu$, $f(xg)=f(x)$ for a.e. $x\in G$. If $G_\mu$ is 
not compact then $\Cal H_{\mu\pi_*}=\{0\}$. 
\endproclaim 
\demo{Proof}When $p>1$, Proposition 2.2 applies. When $p=1$, embed $L^1(G)$ 
in $M(G)$ and use Corollary 2.3. \qed
\enddemo

Another condition that ensures that the Choquet-Deny equation has only 
trivial solutions is that of strict convexity. The proof of the following 
proposition is a routine exercise. 

\proclaim{Proposition 2.5} If the Banach space $E$ (resp., $E_*$) is strictly 
convex then the Choquet-Deny equation in $E$ (resp., $E_*$) admits only 
trivial solutions.  
\endproclaim 
\flushpar 
Thus, in particular, the Choquet-Deny equation in a Hilbert space or in 
an $L^p$-space with $1<p<\iy$ is not interesting. 

We will now introduce two examples 
which provided the original stimulus for our investigations. 

\example{Example 2.6} Let $\Cal X$ be a locally compact $G$-space where $G$ 
is a locally compact group and the mapping $G\times \Cal X\ni (g,x)\to 
gx\in\Cal X$ is continuous. The action of $G$ on $\Cal X$ gives rise to a 
strongly continuous representation $\pi_*$ of $G$ in $C_0(\Cal X)$, defined 
by $\bigl(\pi_*(g)f\bigr)(x)=f(g^{-1}x)$. The adjoint of $\pi_*$ is the 
representation $\pi$ in $M(\Cal X)=C_0(\Cal X)^*$ given by 
$\bigl(\pi(g)\si\bigr)(A)=\si(g^{-1}A)$. 

Suppose $\mu\in M_1(G)$ and $\si\in M_1(\Cal X)$. When both $\Cal X$ and $G$ 
are second countable, $\pi(\mu)\si$ has a simple probabilistic 
interpretation: Let $X$ and $Y$ be independent random variables taking 
values in $\Cal X$ and $G$, resp.. If $\si$ is the distribution of $X$ 
and $\mu$ is the distribution of $Y$ then $\pi(\mu)\si$ is the distibution of 
the $\Cal X$-valued random variable $YX$. The solutions of the Choquet-Deny 
equation $\pi(\mu)\si=\si$ in $M_1(\Cal X)$ are called {\it $\mu$-stationary}
 measures. From the probabilistic point of view they are, precisely, the 
possible limits in $M_1(\Cal X)$, in the weak* topology, of probability 
distributions of 
sequences of $\Cal X$-valued random variables of the form $Y_nY_{n-1}\dots 
Y_1X$ where $\{Y_n\}_{n=1}^\iy$ is a sequence of independent identically 
distributed $G$-valued random variables whose common distribution is $\mu$, 
and $X$ is an $\Cal X$-valued random variable, independent of 
$Y_1,Y_2,\dots$\,. The sequence $Y_nY_{n-1}\dots Y_1X$ forms a Markov chain in 
$\Cal X$ [37, Proposition 4.4, p.\,30], [10, Sec.\,II]. The 
solutions of the Choquet-Deny equation $\pi_*(\tilde\mu)f=f$ in $C_0(\Cal X)$ 
are harmonic functions of this Markov chain. When the Choquet-Deny theorem 
holds, the $\mu$-stationary measures are just the $G_\mu$-invariant 
probability measures on $\Cal X$, and the elements of $\Cal 
H_{\tilde\mu\pi_*}$ are functions $f\in C_0(\Cal X)$ with $f(gx)=f(x)$ for 
every $g\in G_\mu$ and $x\in\Cal X$. \ee
\endexample

\example{Example 2.7} Let $\rho$ be a continuous unitary representation of 
$G$ in a Hilbert space $\frak H$. Then the formula 
$\pi(g)A=\rho(g)A\rho(g^{-1})$ defines a representation $\pi$ of 
$G$ in $B(\frak H)$, the algebra of bounded linear operators on $\frak H$. 
The subrepresentation $\pi_*$ of $\pi$ in the ideal $\Cal T(\frak H)$ of 
trace class operators is strongly continuous (with respect to the trace 
norm) and $\pi$ is the adjoint of $\pi_*$. When the Choquet-Deny theorem 
holds then $\Cal H_{\mu\pi} =\rho(G_\mu)'$, the commutant of $\rho(G_\mu)$ 
in $B(\frak H)$, and $\Cal H_{\mu\pi_*} =\rho(G_\mu)'\cap\Cal T(\frak H)$. 

When $\frak H=L^2(G)$ and $\rho$ is the right regular representation, the 
Choquet-Deny equation in $B(L^2(G))$ can be viewed as a 
``non-commutative'' or ``quantized'' version of the classical Choquet-Deny 
equation. We will refer to $\pi$ as the right regular representation of $G$ 
in $B(L^2(G))$. Of course, $\pi$ commutes with a similarly defined left 
regular representation $\pi_l$ and therefore $\Cal H_{\mu\pi}$ is always 
invariant under $\pi_l$. Moreover, since $B(L^2(G))$ contains a copy of 
$L^\iy(G)$, $\Cal H_{\mu\pi}$ contains a copy of $\Cal H_\mu$. 
When $G_\mu=G$, the subspace of trivial solutions of the Choquet-Deny 
equation in $B(L^2(G))$ is just $VN(G)$, the von Neumann algebra 
generated by $\rho_l$, the left regular representation of $G$ in $L^2(G)$; 
$\Cal H_{\mu\pi}= VN(G)$ if and only if the conclusion of the 
Choquet-Deny theorem is true for $\mu$. The equality $\Cal H_{\mu\pi}= 
VN(G)$ is possible only when $G$ is an amenable group, cf. [26,38]. 

We note that, analogously as in the classical case, the Choquet-Deny equation 
in $\Cal T(L^2(G))$ has always only trivial solutions and when $G_\mu$ is 
not compact, then $\Cal H_{\mu\pi_*}=\{0\}$. This follows from Proposition 
2.2 when $\Cal T(L^2(G))$ is considered as a locally convex space under the 
$\si\bigl(\Cal T(L^2(G)),\Cal K(L^2(G))\bigr)$-topology, where $\Cal 
K(L^2(G))$ denotes the ideal of compact operators. \ee 
\endexample

In connection with Example 2.7 we would like to mention that a very 
different route to establish a ``non-commutative'' 
Choquet-Deny equation has been recently taken 
by Chu and Lau in [5]. There, the duality between $L^\infty(G)$ and $L^1(G)$ 
is replaced by that between the group von Neumann algebra $VN(G)$ and 
the Fourier algebra $A(G)$. Moreover, the measure algebra $M(G)$ is replaced 
by the Fourier-Stieltjes algebra $B(G)$, and probability measures correspond 
to the elements of the set $P^1(G)$ consisting of bounded continuous 
positive definite functions $\si$ on $G$ with $\sigma (e) = 1$. 
The role of the right regular representation of $M(G)$ in $L^\iy(G)$ is 
played by the canonical action of $B(G)$ on $VN(G)$. For 
a fixed $\sigma \in B(G)$, the authors study the space 
$
H_\si=\{ T \in VN(G)\,;\,\sigma T = T \}
$ 
of {\it $\sigma$-harmonic functionals\/} on $A(G)$. Their investigations show 
analogies with but mainly reveal crucial differences from the classical 
situation. For instance, for $\sigma\in P^1(G)$, $H_\sigma$ is always a 
subalgebra of $VN(G)$ -- which is in contrast to the fact that for a 
probability measure $\mu$, the space of 
bounded $\mu$-harmonic functions is a subalgebra of $L^\iy(G)$ only if it 
is trivial. This completely different behaviour is not surprising 
since the classical theory of harmonic functions can be recovered from the 
setting of [5] only for abelian groups. 

\head{\bf 3. The preannihilator of $\Cal H_{\mu\pi}$}\endhead
Given a probability measure $\mu$ on a locally compact group $G$, let $J_\mu$ 
denote the set 
$$
J_\mu=\ov{\{\vp -\vp *\mu\,;\, \vp\in L^1(G)\}}\tag{3.1}
$$
where the bar means closure with respect to the $L^1$-norm. $J_\mu$ is 
evidently a left ideal in the group algebra $L^1(G)$, whose annihilator in 
$L^\iy(G)$ is precisely the space $\Cal H_\mu$ of the bounded $\mu$-harmonic 
functions. As pointed out by Willis [42,43], ideals of this form appear 
naturally not only in the theory of $\mu$-harmonic functions but also in 
the study of amenability and certain factorization questions in group 
algebras. The quotient $L^1(G)/J_\mu$ turns out to be an abstract $L^1$-space 
whose pointwise realization is the boundary needed to represent the 
$\mu$-harmonic functions by means of a Poisson formula. The space $\Cal J$ of 
all ideals of the form $J_\mu$, where $\mu$ ranges over $M_1(G)$, has an 
interesting order structure when ordered by inclusion [42,17]. 
We note that 
$$
J_\mu\sub L^1_0(G,G_\mu)\sub L^1_0(G) 
$$
where $L^1_0(G,G_\mu)$ 
denotes the kernel of the canonical mapping from $L^1(G)$ to \linebreak $L^1(G/G_\mu)$, 
and $L^1_0(G)=L^1_0(G,G)$ is the augmentation ideal $L^1_0(G)=\{f\in 
L^1(G)\,;\,\int f=0\}$. Moreover, $L^1_0(G,G_\mu)$ coincides 
with the preannihilator of the subspace of trivial solutions of the 
Choquet-Deny equation and the equality $J_\mu=L^1_0(G,G_\mu)$ holds if and 
only if the conclusion of the Choquet-Deny theorem is true for $\mu$. 

Let $\pi_*$ be a strongly continuous representation of $G$ by isometries in a 
Banach space $E_*$ and $\pi$ the adjoint of $\pi_*$ acting in the dual $E$ of 
$E_*$. The analog of $J_\mu$ is the closed subspace $J_{\mu\pi}$ of $E_*$ 
given by 
$$
J_{\mu\pi}=\ov{\{x_* -\pi_*(\tilde\mu)x_*\,;\,x_*\in E_*\}} 
\tag{3.2}
$$
where the bar denotes now the norm closure in $E_*$. It is evident that, as 
in the classical case, $\Cal H_{\mu\pi} = J_{\mu\pi}^{\,\perp}$. 

\example{\bf Remark 3.1} The subspace of trivial solutions of the 
Choquet-Deny equation in $E$ is the annihilator of the closed subspace 
$E_{*0}(G_\mu)$ of $E_*$ spanned by the set $\{x_* -\pi_*(g)x_*\,;\,g\in 
G_\mu, x_*\in E_*\}$. $E_{*0}(G_\mu)$ thus plays the role of 
$L^1_0(G,G_\mu)$; $J_{\mu\pi}\sub E_{*0}(G_\mu)$ with equality if and only if 
the Choquet-Deny theorem in $E$ has only trivial solutions. 
\endexample

Throughout the sequel it will be convenient to identify $L^1(G)$ with 
the subspace of absolutely continuous complex measures in $M(G)$. With this 
convention in force we obtain the following simple relation between $J_\mu$ 
and $J_{\mu\pi}$. Recall that for every $\si\in M(G)$, $\tilde\si$ denotes 
the measure  $\tilde\si (A)=\si (A^{-1})$. 

\proclaim{Theorem 3.2} $J_{\mu\pi}=\pi_*(\widetilde{J_\mu})E_* 
=\{\pi_*(\tilde\vp)x_*\,;\,\vp\in J_\mu,\ x_*\in E_*\}$ and 
$J_{\mu\pi}=\pi_*(\widetilde{J_\mu})J_{\mu\pi}$. 
\endproclaim 
\demo{Proof}Given $\vp\in J_\mu$, $x_*\in E_*$, and $x\in \Cal H_{\mu\pi}$, 
we obtain 
$$
\<\pi_*(\tilde\vp)x_*\,,x\>=\<x_*\,,\pi(\vp)x\> =\int_G 
\<x_*,\pi(g)x\>\,\vp(dg) =0, 
$$
because the function $g\to \<x_*,\pi(g)x\>$ is $\mu$-harmonic. Hence, 
$\pi_*(\widetilde{J_\mu})E_*\sub J_{\mu\pi}$. In particular, this means that 
$\pi_*(\widetilde{J_\mu})J_{\mu\pi}\sub J_{\mu\pi}$ and so $J_{\mu\pi}$ is a 
left Banach $\widetilde{J_\mu}$-module. 
But as pointed out in [42, p.\,203], 
if $\{\vep_\al\}_{\al\in A}$ is a bounded approximate identity in $L^1(G)$ 
then 
$$
\eta_{\al n}=\vep_\al *(\del_e-\tfrac1n\sum_{i=1}^n\mu^i),\quad \al\in A,
\ n\in \Bbb N, 
$$
is a bounded right approximate identity for $J_\mu$. Since 
$$
\lim_{n\to\infty}\|\pi_*(\tfrac1n\sum_{i=1}^n\mu^i)x_*\|=0 
$$ 
for every $x_*\in J_{\mu\pi}$, it easily follows that the 
$\widetilde{J_\mu}$-module $J_{\mu\pi}$ has a bounded left approximate 
identity. Hence, Cohen's factorization theorem yields 
$\pi_*(\widetilde{J_\mu})J_{\mu\pi}= J_{\mu\pi}$. Thus 
$J_{\mu\pi}=\pi_*(\widetilde{J_\mu})E_*$. 
\qed\enddemo

\example{Example 2.7 (continued)} Suppose that $E_*$ is a Banach algebra with 
multiplication denoted by $\star\,$. It is immediate that if the identity 
$$
\pi_*(\si)(x_*\star y_*)=x_*\star(\pi_*(\si)y_*)
$$ 
holds for all 
$x_*,y_*\in E_*$ and $\si\in M(G)$, then $J_{\mu\pi}$ will be a left ideal in 
$E_*$. This is, of course, the case when $\pi_*$ is the right regular 
representation in $L^1(G)$ and $\star\,$ is the convolution in 
$L^1(G)$. Continuing our discussion of the ``non-commutative'' 
Choquet-Deny equation in $B(L^2(G))$, we wish to point out here that the 
convolution in $L^1(G)$ has, in fact, an analog in $\Cal T(L^2(G))$ and that 
with respect to this ``non-commutative'' convolution $J_{\mu\pi}$ is a left 
ideal in $\Cal T(L^2(G))$, exactly as in the classical case. Such 
``non-commutative'' convolution was recently introduced and studied in 
[31,32], and we believe that ideals of the form $J_{\mu\pi}$ may be useful in 
further investigation of the resulting ``non-commutative'' version of the 
group algebra $L^1(G)$ [36].

The convolution in $\Cal T(L^2(G))$ can be defined as follows. Let $\pi_l$ 
denote the left regular representation of $G$ in $B(L^2(G))$ and $\pi_{l*}$ 
its restriction to $\Cal T(L^2(G))$. Recall that there is a 
canonical mapping $\k\:\Cal T(L^2(G))\to L^1(G)\sub M(G)$, which commutes 
with the left and right regular representations and is given by $\k (S)(A) 
=\tr[S F(A)]$ for every Borel set $A\sub G$, where $F(A)\in B(L^2(G))$ 
is the operator of multiplication by the characteristic function of $A$; $\k$ 
is the preadjoint of the canonical embedding of $L^\iy(G)$ in $B(L^2(G))$. 
One then defines convolution of two trace class operators $S$ and $T$ by 
$$
S * T = \pi_{l*}\bigl(\k(S)\bigr)T .\tag{3.3}
$$
Since $\k(S * T)=\k(S) *\k(T)$, it easily follows that $*$ is an 
associative algebra product with the further property that 
$$
\tr(S*T) =\tr(S)\tr(T)
$$ 
for all $S,T\in \Cal T(L^2(G))$. $\bigl(\Cal T(L^2(G)),*\bigr)$ is a Banach 
algebra with respect to the trace norm. Since the left and right 
regular representations commute it is also evident that 
$$
\pi_*(\si)(S*T) =S*(\pi_*(\si)T)
$$ 
for all $\si\in M(G)$ and $S,T\in \Cal T(L^2(G))$. Hence, we obtain: 

\proclaim{Proposition 3.3} $J_{\mu\pi}$ is a left ideal in $\bigl(\Cal 
T(L^2(G)),*\bigr)$.\footnotemark
\endproclaim 
\footnotetext{\ The definition of the ``non-commutative'' convolution given 
here is the reverse of the definition originally given in [31,32], i.e., our 
$S*T$ is $T*S$ in [31,32]; this ensures that  
$J_{\mu\pi}$ is, as the classical 
$J_\mu$, a {\it left\/} ideal.}

In connection with Remark 3.1 we note that for $E_*= \Cal T(L^2(G))$, 
$E_{*0}(G)$ is the preannihilator $VN(G)_\perp$ of $VN(G)$ in $\Cal 
T(L^2(G))$ and is properly contained in the augmentation ideal of 
$\bigl(\Cal T(L^2(G)),*\bigr)$, unless $G=\{e\}$. 
 
The ``non-commutative'' convolution can be 
also defined, in essentially the same way, in the space $\Cal N(L^p(G))$ of 
nuclear operators on $L^p(G)$, when $1<p<\iy$. Since $\Cal N(L^p(G))$ is the 
predual of $B(L^p(G))$ it will remain true that the preannihilator of the 
space of solutions of the Choquet-Deny equation in $B\bigl(L^p(G)\bigr)$ is a 
left ideal in $\Cal N(L^p(G))$. \ee 
\endexample 

As we already mentioned, in the classical case the quotient $L^1(G)/J_\mu$ is 
an abstract $L^1$-space whose pointwise realization is the boundary needed to 
represent the $\mu$-harmonic functions by means of a Poisson formula. This 
fact can be established by a purely functional analytic argument, see [42]. 
It is therefore conceivable that a purely functional analytic argument could 
also be used to obtain a Poisson formula for the $\mu$-harmonic vectors and 
relate it to the formula for the classical bounded $\mu$-harmonic 
functions. However, we will not attempt to pursue this approach here, 
choosing instead a more basic, probabilistic approach based on martingale 
theory. The price to be paid for this are the assumptions of separability 
that we will need to impose on the predual $E_*$ of our Banach space and, for 
our main results, also on the group $G$. On the other hand, the probabilistic 
approach yields certain powerful convergence results which do not seem 
possible to obtain by other means. 

\head{\bf 4. Random walks and their harmonic functions}\endhead
In this section we review elements of the classical theory of the bounded 
$\mu$\nolinebreak-\nolinebreak harmonic functions and those elements of the 
theory of random walks 
that are needed to develop our generalization. With the exception of some 
more specialized results pertaining to the theory of $\mu$-boundaries, 
most of the material presented here is a special case of the basic theory of 
Markov chains to be found, e.g., in [33] or [37]. 

Let $G$ be a locally compact group and $\mu$ a probability measure on $G$. By 
the (right) random walk of law $\mu$ one means the Markov chain with state 
space $G$ and transition probability $\Pi(g,A)=\mu(g^{-1}A)$. The position of 
the random walk after its $n$-th step ($n=0,1,\dots$) can be expressed as the 
product $Y_0Y_1\dots Y_n$ where $Y_0$ is the initial position and 
$Y_1,Y_2,\dots$ are independent $G$-valued random variables distributed 
according to the law $\mu$. In general, $Y_0$ is also a random variable, 
independent of $Y_1,Y_2,\dots$ and distributed according to a law $\nu$. 

Let $G^\iy$ denote the product space $G^\iy =\prod_{n=0}^\infty G$ (the space 
of paths $\om=\{\om_n\}_{n=0}^\infty$ of the random walk), 
$X_n\:G^\infty\to G$, $n=0,1,\dots$\,, the canonical projections, 
and $\Cal B^\infty$ the product $\si$-algebra $\prod_{n=0}^\iy \Cal 
B\alb=\si (X_0,X_1,\dots)$, where $\Cal B$ stands for the $\si$-algebra of 
Borel subsets of $G$. The law $\mu$ of the random walk and the starting 
measure $\nu$ define a measure $Q_\nu$ on $(G^\iy,\Cal B^\iy)$, called the 
{\it Markov measure\/}. The canonical projections $X_n$ become random 
variables on the probability space $(G^\iy,\Cal B^\iy,Q_\nu)$, with 
distributions $\nu*\mu^n$. When $G$ is second countable, $Q_\nu$ is the image 
of the product measure $\nu\ti\mu\ti\mu\ti\dots$ on $(G^\iy,\Cal B^\iy)$ 
under the mapping $G^\iy\ni (\om_0,\om_1,\om_2,\dots)\to 
(\om_0,\om_0\om_1,\om_0\om_1\om_2,\dots)\in G^\iy$. 
\footnote{\ Second countability ensures that the mapping 
$(\om_0,\om_1,\om_2,\dots)\to (\om_0,\om_0\om_1,\om_0\om_1\om_2,\dots)$ is 
measurable. To ensure this without second countability requires extending the 
infinite product measure $\nu\ti\mu\ti\mu\ti\dots$ to a $\sigma$-algebra 
larger than $\Cal B^\iy$.}

In the case that $\nu$ is a point measure $\del_g$, we will write $Q_g$ 
rather than $Q_{\del_g}$. The function $G\ti \Cal B^\iy\ni (g,A)\to Q_g(A)$ 
is, in fact, a transition probability from $(G,\Cal B)$ to $(G^\iy,\Cal 
B^\iy)$, and one has 
$$
Q_\nu(A)=\int_G Q_g(A)\,\nu(dg),\quad A\in\Cal B^\iy. \tag{4.1}
$$

The transition probability $\Pi(g,A)=\mu(g^{-1}A)$ is invariant with respect 
to the action of $G$ on $G$ by left translations, i.e., $\Pi(gg',gA) 
=\Pi(g',A)$ for all $g,g'\in G$ and $A\in\Cal B$. There is also a related 
action of $G$ on the space of paths $G^\iy$, namely, $g\{\om_n\}_{n=0}^\iy 
=\{g\om_n\}_{n=0}^\iy$. With this action $(G^\iy,\Cal B^\iy)$ is a Borel 
$G$-space. The Markov measures $Q_g$ satisfy 
$$
Q_{gg'}(gA)=Q_{g'}(A),\quad g,g'\in G,\ A\in\Cal B^\iy\tag{4.2}
$$
(i.e., the transition probability $(g,A)\to Q_g(A)$ is $G$-invariant). 

In general, when $f\:\Cal X\to\Cal X'$ is a Borel function from a Borel 
space $(\Cal X,\Cal A)$ into a Borel space $(\Cal X',\Cal A')$ and $\si$ a 
measure on $(\Cal X,\Cal A)$, we will write $f\si$ for the measure 
$(f\si)(A')=\si(f^{-1}(A'))$, $A'\in \Cal A'$. When $(\Cal X,\Cal A)$ is a 
Borel $G$-space and $g\in G$, then $g\si$ will stand for the measure 
$(g\si)(A)=\si(g^{-1}A)$, $A\in \Cal A$. The convolution $\nu *\si$ of a 
measure $\nu$ on $G$ with $\si$ is defined by 
$$
(\nu *\si)(A)=\int_G(g\si)(A)\,\nu(dg),\quad A\in\Cal A\tag{4.3}
$$
(provided that the function $G\ni g\to (g\si)(A)$ is Borel for every 
$A\in\Cal A$). Using the above notation we have $Q_g=gQ_e$ (Eq.\,(4.2)) and 
$Q_\nu= \nu *Q_e$ (Eq.\,(4.1)). 

Let $\Bbb B(G)$ and $\Bbb B(G^\iy)$ denote the algebras of bounded 
complex-valued Borel functions on $G$ and $G^\iy$, resp., equipped 
with the sup-norms. The natural action of $G$ on $\Bbb B(G)$ and $\Bbb 
B(G^\iy)$ is 
the action $(gf)(x)=f(g^{-1}x)$, associated with the left regular 
representations. It follows from Eq.\,(4.2) that the formula 
$$
(Rf)(g)=\int_{G^\iy}f(\om)\,Q_g(d\om)=\int_{G^\iy}f(g\om)\,Q_e(d\om)\tag{4.4}
$$
defines an equivariant contraction $R\:\Bbb B(G^\iy)\to\Bbb B(G)$. 

We will denote by $r$ the right regular representation of $G$ in $\Bbb B(G)$, 
as well as its extension to a representation of $M(G)$. 

The {\it Markov shift\/} $\vt\,$ is the transformation $\vt\:G^\iy\to G^\iy$ 
given by $\vt(\{\om_n\}_{n=0}^\iy) =\{\om_{n+1}\}_{n=0}^\iy$. $\vt$ commutes 
with the natural $G$-action on $G^\iy$ and transforms the Markov measure 
$Q_\nu$ into the Markov measure $Q_{\nu*\mu}$, i.e.,
$$
\vt Q_\nu=Q_{\nu*\mu}.\tag{4.5}
$$
The formula $\theta f=f\circ\vt$ defines an injective homomorphism $\theta$ 
of $\Bbb B(G^\iy)$ into itself. It is an immediate consequence of Eqs (4.1) 
and (4.5) that 
$$
R\theta=r(\mu)R. \tag{4.6}
$$

The $\si$-algebra $\Cal B^{(i)}=\{A\in\Cal B^\iy\,;\,\vt^{-1}(A)=A\}$ is 
called the {\it invariant $\si$-algebra}. Elements of $\Cal B^{(i)}$ are 
called {\it invariant sets\/}, and $\Cal B^{(i)}$-measurable functions are 
called 
{\it invariant random variables\/}. A\phantom, $\Cal B^\iy$-measurable 
function $f\:G^\iy \to\Bbb C$ is $\Cal B^{(i)}$-measurable if and only if 
$f\circ \vt= f$. Since $\vt$ and the $G$-action on $G^\iy$ commute, $\Cal 
B^{(i)}$ is preserved by the $G$-action and therefore $(G^\iy,\Cal B^{(i)})$ 
is also a Borel $G$-space. 

We will denote by $\Bbb B_i$ the algebra of bounded 
complex-valued invariant random variables. $\Bbb B_i$ is precisely the 
subalgebra of $\Bbb B(G^\iy)$ consisting of the fixed points of the 
homomorphism $\theta\:\Bbb 
B(G^\iy)\to\Bbb B(G^\iy)$. By Eq.\,(4.6), for every $f\in 
\Bbb B_i$, $r(\mu)Rf=Rf$, i.e., $Rf$ is a bounded $\mu$-harmonic function. We 
note that this property of $R$ is equivalent to having 
$$
Q_e(A)=(\mu* Q_e)(A)\tag{4.7}
$$
for every $A\in \Cal B^{(i)}$. 

We will denote by $\Bbb H_\mu$ the space of 
bounded $\mu$-harmonic functions in $\Bbb B(G)$, equipped with the 
sup norm. $\Bbb H_\mu$ is invariant under the left regular representation; 
by the action of $G$ on $\Bbb H_\mu$ we will always mean the action 
associated with the left regular representation. 

An invariant set $A\in\Cal B^{(i)}$ will be called {\it universally null\/} 
(resp., {\it universally conull\/}) if $Q_g(A)=0$ (resp., $Q_g(A)=1$) 
for every $g\in G$. We will say that a property dependent on $\om\in G^\iy$ 
holds {\it universally almost everywhere\/} (u.a.e.) if it holds for $\om$ in 
a universally conull set.

Let $\Cal N_u$ denote the collection of universally null sets. An 
invariant random variable $f\:G^\infty\to\Bbb C$ will be called {\it 
universally essentially bounded\/} if 
$$
\|f\|_u =\inf_{\Del\in\Cal N_u}\bigl(\,\sup_{\omega\in G^\infty 
-\Del}|f(\om)|\bigr)<\iy .\tag{4.8}
$$
$\|\cdot\|_u$ is a $C^*$-norm on the *-algebra $\Bbb L^\iy_i(\mu)$ of 
equivalence classes of the universally essentially bounded invariant random 
variables, where two such random variables are equivalent when they coincide 
u.a.e.. Since $\Cal N_u$ is invariant under the action of $G$ on $G^\iy$, the 
natural action of $G$ on $\Bbb B_i$ (the left regular representation) induces 
an action of $G$ on $\Bbb L^\iy_i(\mu)$, and the contraction $R$ of 
Eq.\,(4.4) induces an equivariant contraction, which we denote also by $R$, 
of $\Bbb L^\iy_i(\mu)$ into $\Bbb H_\mu$. The following fundamental result 
[33, Proposition V.2.4] is a well known consequence of the Martingale 
Convergence Theorem.

\proclaim{Proposition 4.1} $R$ is an equivariant isometric isomorphism of 
\phantom{!}$\Bbb L^\iy_i(\mu)$ onto $\,\Bbb H_\mu$. Moreover, for every 
$h\in\Bbb H_\mu$ 
the sequence $\{h\circ X_n\}_{n=0}^\infty$ converges u.a.e. to $R^{-1}h$.
\endproclaim 

We will find it useful in our treatment of the $\mu$-harmonic vectors in 
Section 6 that the sequence $\{h\circ X_n\}_{n=0}^\infty$ converges u.a.e. 
not only when $h\in \Bbb H_\mu$ but also when $h$ is a 
$\mu$\nolinebreak-\nolinebreak subharmonic function. A bounded above Borel 
function $h\: G\to \Bbb R$ is called {\it 
$\mu$-subharmonic\/} if 
$$
h(g)\le\int_G h(gg')\,\mu(dg')\tag{4.9}
$$ 
holds for every $g\in G$. Clearly, if $h$ is bounded harmonic then $|h|$ is 
subharmonic. It is also easy to see that if $h_1,\dots,h_n$ are subharmonic 
functions then so is $h=\max_{1\le i\le n}h_i$. The next proposition is a 
direct consequence of the theory of submartingales.

\proclaim{Proposition 4.2} Given a bounded above subharmonic function 
$h\:G\to \Bbb R$ the sequence $h\circ X_n$ converges u.a.e. to an invariant 
random variable $f\:G^\iy\to \Bbb R$ such that for every $g\in G$, 
$$
h(g)\le\int_{G^\iy}f(\om)\,Q_g(d\om).
$$
\endproclaim 
\phantom{bbbbbbbbbbbbbbbB}
                                  
Returning to the $\mu$-harmonic functions, observe that since 
$\Bbb L^\iy_i(\mu)$ is an abelian $C^*$-algebra, $\Bbb H_\mu$ itself is an 
abelian $C^*$-algebra when equipped with the product 
$$
h_1\di h_2 
=R[(R^{-1}h_1)(R^{-1}h_2)].
$$ 
The following intrinsic description of this product is an immediate 
consequence of Proposition 4.1 and the Dominated Convergence Theorem.

\proclaim{Corollary 4.3} Given $h_1,h_2\in \Bbb H_\mu$, the sequence 
$r(\mu^n)(h_1h_2)$ converges pointwise to $h_1\diamond h_2$. 
\endproclaim 
\flushpar
In general, $\di$ differs from the usual pointwise product of functions. 
It is not hard to see that $\di$ coincides with the usual product if and only 
if for every $A\in \Cal B^{(i)}$ and every $g\in G$, $Q_g(A)$ is either 0 or 
1 (since $Q_g=gQ_e$ this is equivalent to having $Q_e(A)\in\{0,1\}$ for every 
$A\in\Cal B^{(i)}$). The Hewitt-Savage 0-1 law [13], [29, p.\,190] implies 
that for second countable abelian groups $\di$ is always identical with the 
usual product. This, in turn, implies the Choquet-Deny theorem, see [29, 
pp.\,192-193]. 

The measure $\mu$ is called {\it spread out\/} if for some $n$ the 
convolution power $\mu^n$ is nonsingular (of course, on a discrete group 
every measure is spread out). Bounded $\mu$-harmonic functions of a spread 
out measure are right uniformly continuous [4, Proposition I.6, p.\,23]. It 
can be shown that for spread out $\mu$, $\di$ is the usual product if and 
only if every bounded $\mu$-harmonic 
function is constant on the left cosets of $G_\mu$, i.e., the Choquet-Deny 
equation in $\Bbb B(G)$ has only trivial solutions. When $G$ is nilpotent, or 
is compactly generated and has polynomial growth, then for every spread out 
$\mu$ the Choquet-Deny equation has only trivial solutions [4, Proposition 
IV.10, p.\,98], [21]. The largest class of countable groups known today for 
which this is the case is the class of FC-hypercentral groups [16, Theorem 
4.8]. On the 
other hand, $\di$ differs from the usual product for every probability 
measure (spread out or not) for which $G_\mu$ is nonamenable.
\footnote{\ When $G_\mu$ is nonamenable, then by, e.g., [10, p.\,213] 
there exists a nonconstant bounded continuous $\mu$-harmonic function 
$h\:G_\mu\to\Bbb C$. A version of the argument in [29, pp.\,192-193] shows 
that then $0<Q_e(A)<1$ for some $A\in \Cal B^{(i)}$. So $\di$ differs from 
the usual product.} 
We also note that every $\si$-compact amenable locally compact group admits 
an absolutely continuous probability measure $\mu$ such that $G=G_\mu$ and 
the Choquet-Deny equation has only trivial solutions [26,38]; but amenability 
of $G$ alone does not guarantee that every absolutely continuous probability 
measure has this property [26], [20, Theorem 3.13], [21, Theorem 3.16]. 

Since every $G$-invariant function $f\:G\to \Bbb C$ is constant, Proposition 
4.1 has also this immediate corollary. 

\proclaim{Corollary 4.4} $G$ acts ergodically on $\Bbb L^\iy_i(\mu)$, 
i.e., if $f\in \Bbb L^\iy_i(\mu)$ and $gf=f$ for every $g\in G$, then $f$ 
is constant u.a.e..
\endproclaim 
             
We will now discuss the Choquet-Deny equation in $L^\iy(G)$. Recall that the 
space of $\mu$-harmonic functions in $L^\iy(G)$ is denoted by $\Cal H_\mu$. 
Evidently, every $h\in \Bbb H_\mu$ defines an element of $\Cal H_\mu$. When 
the Haar measure $\lam$ is $\si$-finite, i.e., when $G$ is $\si$\nolinebreak 
-\nolinebreak compact, one can prove using martingales that every element of 
$\Cal H_\mu$ arises in this way, i.e., $\Cal H_\mu$ is precisely the space of 
equivalence classes of the elements of $\Bbb H_\mu$ modulo $\lam$ 
[9, Proposition 2]. 
$\si$-finiteness of $\lam$ is also used in the proofs of Proposition 4.5 and 
Corollary 4.6 below. {\bf We will therefore assume for the remainder of this 
section that $G$ is $\si$-compact.} 

For spread out $\mu$ every $h\in\Bbb H_\mu$ is continuous. Therefore for 
such $\mu$, $\Cal H_\mu=\Bbb H_\mu$. However, in general, $\Cal H_\mu$ can be 
very different from $\Bbb H_\mu$. For example, when $G$ is abelian and $\mu$ 
is any discrete probability measure with $G_\mu=G$, then $\Cal H_\mu=\Bbb C1$ 
by the Choquet-Deny theorem; at the same time every bounded Borel function 
constant on the cosets of the subgroup generated by the discrete support of 
$\mu$ is $\mu$-harmonic. 

By Eqs\,(4.1) and (4.2) the Markov measure $Q_\lam$, where $\lam$ is the left 
Haar measure, is a $G$-invariant measure 
on $(G^\iy,\Cal B^\iy)$ and therefore the natural action of $G$ on $G^\iy$ 
induces an action of $G$ on $L^\infty(G^\iy,\Cal B^\iy,Q_\lam)$ and on 
$L^\iy(G^\iy,\Cal B^{(i)},Q_\lam)$. We will denote the latter space by 
$L^\iy_i(\mu)$. Eq.\,(4.4) defines an equivariant 
contraction which (abusing notation) we will still denote by $R$, of 
$L^\iy_i(\mu)$ into $L^\iy(G)$. It is clear that $RL^\iy_i(\mu)\sub\Cal 
H_\mu$. The following analog of Proposition 4.1 is also a consequence of the 
theory of martingales. 

\proclaim{Proposition 4.5} $R$ is an equivariant isometric isomorphism 
of $L^\iy_i(\mu)$ onto $\Cal H_\mu$. For every $h\in \Cal H_\mu$ the sequence 
$\{h\circ X_n\}_{n=0}^\iy$ converges $Q_\lam$-a.e. to $R^{-1}h$.
\endproclaim 

We note that, in general, the measure $Q_\lam$ fails to be $\si$-finite when 
restricted to $\Cal B^{(i)}$.
\footnote{\ When $\mu$ is spread out, $\si$-finiteness of $Q_\lam$ implies 
that every bounded $\mu$-harmonic function is constant on the cosets of 
$G_\mu$, see [20, Proposition 2.6] and [21, Lemma 2.3].} 
However, by Eq.\,(4.1) the measure class of the Markov measure $Q_\nu$ is 
completely determined by the measure class of $\nu$. Thus when $\be$ is any 
finite measure equivalent to $\lam$, then $Q_\be$ is a finite measure 
equivalent to $Q_\lam$ and one can replace $Q_\lam$ by $Q_\be$. Consequently, 
$L^\iy_i(\mu)$ and $\Cal H_\mu$ are abelian $W^*$-algebras. Denoting by $\di$ 
the product 
$$
h_1\di h_2= R[(R^{-1}h_1)(R^{-1}h_2)]
$$ 
in $\Cal H_\mu$, one can easily prove the following analog of Corollary 
4.3. $r$ denotes now the right regular representation in $L^\iy(G)$. 

\proclaim{Corollary 4.6} Given $h_1,h_2\in\Cal H_\mu$, the sequence 
$r(\mu^n)(h_1h_2)$ converges to $h_1\di h_2$ $\lam$-a.e. and (hence) in the 
weak* topology. 
\endproclaim 
\example{\bf Remark 4.7} Let $L^1_i(\mu)$ denote the space of complex 
measures on $(G^\iy,\Cal B^{(i)})$, absolutely continuous with respect to 
$Q_\lam$. Note that if $\vp\in L^1(G)\sub M(G)$ then $\vp *Q_e\in L^1_i(\mu)$ 
and $\<\vp\,,Rf\>= \<\vp *Q_e\,,f\>$ for every $f\in L^\iy_i(\mu)$. Hence, 
$R$ is weakly* continuous as a mapping into $L^\iy(G)$ and the preadjoint 
$R_*\: L^1(G)\to L^1_i(\mu)$ of $R$ is given by $R_*\vp=\vp*Q_e$. 
In fact, by Lemma 4.8 below, 
the inverse mapping $R^{-1}\:\Cal H_\mu \to L^\iy_i(\mu)$ is also 
weakly* continuous. Thus the weak* (ultraweak) topology of $\Cal H_\mu$ 
considered as the $W^*$-algebra coincides with the restriction of the weak* 
topology of $L^\iy(G)$ to $\Cal H_\mu$. 
The kernel of $R_*$ is precisely the preannihilator $J_\mu$ of $\Cal H_\mu$ 
in $L^1(G)$. It follows that $L^1(G)/J_\mu$ is isometrically isomorphic to 
$L^1_i(\mu)$, i.e., $L^1(G)/J_\mu$ is an abstract $L^1$-space, a result 
already mentioned in Section 3. 
\endexample 

\proclaim{Lemma 4.8} Let $X$ and $\,Y$ be Banach spaces and $T$ a weakly* 
continuous isometry of a weakly* closed subspace $V\sub X^*$ into $Y^*$. Then 
$W=TV$ is also weakly* closed and $T^{-1}\:W\to V$ is weakly* continuous. 
\endproclaim 

\demo{Proof}The proof is a routine application of the Krein-Smulian Theorem 
(see, e.g., [8, Theorem 7, Chap.\,V.5 and Corollary 11, Chap.\,V.3]). 
\qed\enddemo
\proclaim{Proposition 4.9} The following conditions are equivalent: 
\roster
\item"{(i)}" The product $\di$ in $\Cal H_\mu$ coincides with the usual 
product in $L^\iy(G)$. 
\item"{(ii)}" For each $A\in \Cal B^{(i)}$, $Q_g(A)\in\{0,1\}$ for 
$\lam$-a.e. $g\in G$.
\item"{(iii)}" The Choquet-Deny theorem is true for $\mu$, i.e., every 
bounded continuous $\mu$-harmonic function is constant on the 
left cosets of $G_\mu$. 
\item"{(iv)}" The Choquet-Deny equation in $L^\iy(G)$ has only trivial 
solutions.
\endroster\endproclaim 
\demo{Proof}That (i)$\LRa$(ii) and (iii)$\LRa$(iv) is straightforward. 

(iv)$\Ra$(i): Let $h_1,h_2\in \Cal H_\mu$. Then for every $\vp\in L^1(G)$, 
$$
\<\vp\,,\pi(\mu^n)(h_1h_2)\>=\int_G\<\vp\,,\pi(g)(h_1h_2)\>\,\mu^n(dg)=
\int_G\<\vp\,,h_1h_2\>\,\mu^n(dg)=\<\vp\,,h_1h_2\>.
$$
Hence, $\pi(\mu^n)(h_1h_2)=h_1h_2$ and thus $h_1\di h_2=h_1h_2$.

(i)$\Ra$(iii): Let $h$ be a bounded continuous $\mu$-harmonic function. Then 
$|h|^2$ is also $\mu$-harmonic. So for every $g\in G$, 
$$
\int_G |h(gt)-h(g)|^2\,\mu(dt)=\int_G |h(gt)|^2-h(gt)^*h(g)-
h(gt)h(g)^*+|h(g)|^2\,\,\mu(dt)=0.
$$
Using continuity of $h$ this implies that $h(gt)=h(g)$ for every 
$t\in\supp\mu$. Thus $r(t)h=h$ for every $t\in\supp\mu$ and, hence, also 
for every $t\in G_\mu$. Therefore $h$ must indeed be constant on the left 
cosets of $G_\mu$. 
\qed\enddemo

Since every universally null set is $Q_\lam$-null, $L^\iy_i(\mu)$ is 
canonically the quotient of $\Bbb L^\iy_i(\mu)$. Continuity of the bounded 
$\mu$-harmonic functions of a spread out measure and Eq.\,(4.1) imply that 
for spread out $\mu$ each of the Markov measures $Q_g$, $g\in G$, is 
absolutely continuous with respect to $Q_\lam$. Consequently, an invariant 
set $A\in\Cal B^{(i)}$ is universally null if and only if it is 
$Q_\lam$-null, and so $\Bbb L^\iy_i(\mu)$ and $L^\iy_i(\mu)$ coincide for 
spread out $\mu$. 
As our discussion of the relation between $\Bbb H_\mu$ and $\Cal H_\mu$ 
indicates, for non-spread out $\mu$, $\Bbb L^\iy_i(\mu)$ and $L^\iy_i(\mu)$ 
can be very different. Moreover, in general, $\Bbb L^\iy_i(\mu)$ is not a  
$W^*$-algebra.
\footnote{\ E.g., when $\mu=\del_e$, then $\Bbb 
L^\iy_i(\mu)\cong \Bbb B(G)$; when $G=\Bbb R$ and $\mu$ is a discrete 
probability measure whose discrete support generates $\Bbb Q$, then $\Bbb 
L^\iy_i(\mu)\cong \Bbb B(\Bbb R/\Bbb Q)$.} 

In contrast to the case of spread out probability measures (or probability 
measures on discrete groups) the Choquet-Deny theorem for general measures 
on continuous noncompact nonabelian locally compact groups has been verified 
only for 2-step nilpotent groups [12], nilpotent [SIN] groups [18], and some 
more special groups [12]. 

The formula 
$$
h(g)=\int_{G^\iy}f(\om)\,Q_g(d\om)
$$ 
is the Poisson formula for 
the bounded $\mu$-harmonic functions and the Borel space $(G^\iy,\Cal 
B^{(i)})$ can be regarded as a ``boundary''. However, this boundary is not a 
``nice'' space , its pathological feature being the fact that 
the $\si$-algebra $\Cal B^{(i)}$ does not separate points (unless $G=\{e\}$). 
Moreover, $(G^\iy,\Cal B^{(i)})$ itself does not depend on $\mu$, this 
dependence being encoded entirely in the properties of the Poisson 
kernel $Q_g$. Apart from the problem of extending the Choquet-Deny theorem 
to nonabelian groups, a major problem in the theory of bounded 
$\mu$-harmonic functions has been to find a realization of the boundary as a 
topological (at least) Hausdorff space whose points can be regarded as limit 
points of the trajectories of the random walk (see [23-25,27] and references 
therein). This theory is outside the scope of this article and we will 
confine ourselves to mentioning only certain generalities regarding the 
possibility to realize the boundary as a ``nice'' space, related to the fact 
that $\Cal H_\mu$ is an abelian $W^*$-algebra. 

By a $G$-space $(\Cal X,\Cal A,\al)$ we mean a Borel $G$-space $(\Cal 
X,\Cal A)$ with a quasiinvariant measure $\al$. 
When $(\Cal X,\Cal A,\al)$ and $(\Cal X',\Cal A',\al')$ are two such 
$G$-spaces we say that $L^\infty(\Cal X,\Cal A,\al)$ is {\it 
$G$-isomorphic\/} to 
$L^\infty(\Cal X',\Cal A',\al')$ if there exists an equivariant *-isomorphism 
of the *-algebra $L^\iy(\Cal X,\Cal A,\al)$ onto $L^\iy(\Cal X',\Cal 
A',\al')$. We define the {\it $\mu$-boundary\/} as any $G$-space $(\Cal 
X,\Cal A,\al)$ such that $L^\infty(\Cal X,\Cal A,\al)$ is $G$-isomorphic to 
$L^\iy_i(\mu)$. We remark that since $Q_\lam$ is equivalent to a finite 
measure, the same is true about the quasiinvariant measure $\al$ on any 
$\mu$-boundary. Moreover, a $G$-space $(\Cal X,\Cal A,\al)$ is a 
$\mu$-boundary if and only if there exists an equivariant identity preserving 
isometry of $L^\iy(\Cal X,\Cal A,\al)$ onto $\Cal H_\mu$. 

Let us assume that $G$ is second countable. Then $L^1(G)$ is separable and 
therefore so is the predual of $\Cal H_\mu$. Hence, $L^1_i(\mu)$ is 
also separable. Moreover, given $f\in L^\iy_i(\mu)$ and $\vp\in L^1_i(\mu)$, 
the function $G\ni g\to \<\vp\,,gf\>\in\Bbb C$ is Borel (in fact, 
continuous). These properties permit a routine application of the classical 
Mackey's theorem about pointwise realization of group actions [28]. 
It shows that there always exist $\mu$-boundaries that are 
standard Borel spaces such that the map $G\ti \Cal X\ni (g,x)\to gx\in 
\Cal X$ is Borel. Moreover, by a theorem of Varadarajan [41, Theorem 3.2] 
one can even take $\Cal X$ to be a compact metric space with the 
map $G\ti\Cal X\ni (g,x)\to gx\in\Cal X$ continuous. If $(\Cal X,\Cal 
A,\al)$ is a standard $\mu$-boundary with the map $G\ti \Cal X\ni (g,x)\to 
gx\in \Cal X$ Borel and $\Phi\:L^\iy(\Cal X,\Cal A,\al)\to L^\iy_i(\mu)$ is 
the equivariant isomorphism, then there exists a probability measure $\rho$ 
on $\Cal X$ (the Poisson kernel) such that 
$$
(R\Phi f)(g)=\int_G f(gx)\,\rho(dx)\pmod\lam
$$ 
for every $f\in L^\iy(\Cal X,\Cal A,\al)$.

Another realization of the $\mu$-boundary, often called the Poisson 
space [4], is the spectrum of the C*-subalgebra of $\Cal H_\mu$ consisting 
of left uniformly continuous bounded $\mu$-harmonic functions; the 
disadvantage of the Poisson space is that it usually is not metrizable.  

Finally, we note that the action of $G$ on a $\mu$-boundary has certain 
very distinctive properties. It is always aproximately transitive and 
amenable [7,14,15,44], and when $\mu$ is spread out it is strongly 
approximately transitive [21]. There are also further properties of a more 
probabilistic nature related to the convergence properties of the underlying 
random walk [10,11,19]. Amenability of the $G$-action on the $\mu$-boundary 
is intimately related to the existence of a norm 1 equivariant projection of 
$L^\iy(G)$ onto $\Cal H_\mu$; in the sequel we will obtain a generalized 
version of this result in the setting of the Choquet-Deny equation in a 
dual Banach space. 

\head{\bf 5. Vector valued harmonic functions}\endhead

Let $E$ be a Banach space with a separable predual $E_*$. 
As is well known, the Borel structure in $E_*$ defined by the norm topology 
is the same as that defined by the weak topology and also the same as 
the weak Borel structure generated by the functions $\<\cdot\,,x\>$, $x\in 
E$. Thus a function $f$ from a Borel space $(\Cal X,\Cal A)$ to $E_*$ 
is measurable (Borel) if and only if for each $x\in E$ the function 
$\Cal X\ni t\to\<f(t)\,,x\>\in \Bbb C$ is Borel. 
We note that if 
$f\:\Cal X\to E_*$ is Borel then the function $t\to \|f(t)\|$ is also Borel. 

Given a complex measure $\si$ on $(\Cal X,\Cal A)$ and a Borel function 
$f\:\Cal X\to E_*$ such that $\int_{\Cal X} \|f\|\,d|\si| <\iy$, by the 
integral $\int_{\Cal X} f\,d\si$ we will mean the unique vector $\int_{\Cal 
X} f\,d\si\in E_*$ with
$$
\bigl<\int_{\Cal X}f\,d\si\,,x\bigr>=
\int_{\Cal X}\<f(t)\,,x\>\,\si(dt)\tag{5.1}
$$ 
for every $x\in E$. 

The proper Borel structure in $E$ to work with is that given by the weak* 
topology. This Borel structure is well known to be standard and coincides 
with the weak Borel structure generated by the functions $\<x_*\,,\cdot\>$, 
$x_*\in E_*$. A function $f$ from a Borel space $(\Cal X,\Cal A)$ to $E$ is 
therefore Borel if and only if for each $x_*\in E_*$ the function $\Cal X\ni 
t\to\<x_*\,,f(t)\>\in \Bbb C$ is Borel. If $f\:\Cal X\to E$ is Borel then the 
function $t\to \|f(t)\|$ is also Borel. 

Given a complex measure $\si$ on 
$(\Cal X,\Cal A)$ and a Borel function $f\:\Cal X\to E$ such that $\int_{\Cal 
X} \|f\|\,d|\si| <\iy$, the integral $\int_{\Cal X} f\,d\si$ is the unique 
vector $\int_{\Cal X} f\,d\si\in E$ with
$$
\bigl<x_*\,,\int_{\Cal X}f\,d\si\bigr>=
\int_{\Cal X}\<x_*\,,f(t)\>\,\si(dt)\tag{5.2}
$$ 
for every $x_*\in E_*$. 

Let $G$ be a locally compact group and $\mu\in M_1(G)$. 
A bounded (with respect to the norm on $E$) Borel function $h\:G\to E$ will 
be called {\it $\mu$-harmonic} if it satisfies Eq.\,(1.1) for every $g\in G$. 
The goal of this section is to extend Proposition 4.1 and Corollary 4.3 to 
such $\mu$-harmonic functions. This will prove very useful in the 
following sections, in our study of the $\mu$-harmonic 
vectors arising from a representation of $G$ by isometries in $E$. 

Let $\Bbb B(G,E)$ and $\Bbb B_i(E)$ denote the Banach spaces of bounded Borel 
functions $f\:G\to E$ and bounded $\Cal B^{(i)}$-measurable functions 
$f\:G^\iy\to E$ ($E$-valued invariant random variables), resp., 
equipped with the sup-norms $\|\cdot\|_{\sup}$. The concept of a universally 
essentially bounded invariant random variable with values in $E$ is 
introduced as in the classical case and the Banach space of equivalence 
classes of such random variables (with norm $\|\cdot\|_u$) will be denoted by 
$\Bbb L^\iy_i(\mu,E)$. Finally, let $r_E$ stand for the right regular 
representation of $G$ in $\Bbb B(G,E)$, as well as its extension to a 
representation of $M(G)$, and let $\Bbb H_\mu(E)$ denote the space of 
bounded $E$-valued $\mu$-harmonic functions. The group acts on $\Bbb B(G,E)$, 
$\Bbb B_i(E)$, $\Bbb L^\iy_i(\mu,E)$, and $\Bbb H_\mu(E)$ in the same way as 
in the classical case. 

It is clear that given $f\in\Bbb B_i(E)$, the function 
$$
h(g)= \int_{G^\iy} f(\om)\,Q_g(d\om)\tag{5.3}
$$
is bounded $\mu$-harmonic, and Eq.\,(5.3) defines an equivariant contraction 
$R_E$ from $\Bbb B_i(E)$ into $\Bbb H_\mu(E)$. As in the classical case we 
will use the same symbol $R_E$ to denote the contraction from 
$\Bbb L^\iy_i(\mu,E)$ into $\Bbb H_\mu(E)$ defined by Eq.\,(5.3). 

\proclaim{Lemma 5.1} If $f\:G^\iy\to E$ is a universally essentially 
bounded invariant random variable and $h=R_Ef$, 
then the sequence $h\circ X_n$ converges weakly* u.a.e. to $f$ and the 
sequence $\|h\circ X_n\|$ converges u.a.e. to $\|f\|$. 
\endproclaim 
\demo{Proof} Given $x_*\in E_*$ consider the functions $h_{x_*}\:G\to \Bbb C$ 
and $f_{x_*}\:G^\iy\to \Bbb C$ defined by $h_{x_*}(g)=\<x_*\,,h(g)\>$ and 
$f_{x_*}(\om)=\<x_*\,,f(\om)\>$, resp.. Then $h_{x_*}\in\Bbb H_\mu$ 
and $h_{x_*}=Rf_{x_*}$. By Proposition 4.1, 
$f_{x_*}=\lim_{n\to\infty}h_{x_*}\circ X_n$ u.a.e.. Since $E_*$ is 
separable this easily implies that $\w*lim_{n\to\iy} h\circ X_n=f$ u.a.e..

The proof of the second statement invokes Proposition 4.2. Let 
$H(g)=\|h(g)\|$, $g\in G$. Then $H$ is a bounded $\mu$-subharmonic function 
and therefore there exists an invariant random variable $F\:G^\iy\to[0,\iy)$ 
with $F=\lim_{n\to\infty}H\circ X_n$ u.a.e.. We need to prove that 
$F(\om)=\|f(\om)\|$ u.a.e.. It is easy to see that $\|f(\om)\|\le F(\om)$ 
u.a.e.. Therefore to prove that $F(\om)=\|f(\om)\|$ u.a.e., it suffices to 
prove that 
$$
\int_{G^\iy}F(\om)\,Q_g(d\om)= \int_{G^\iy}\|f(\om)\|\,Q_g(d\om)
$$ 
for every $g\in G$. 

Let $\{x_{*i}\}_{i=1}^\iy$ be a sequence dense in the unit ball of $E_*$. 
For each $k=1,2,\dots$\,, define a function $H_k\: G\to [0,\iy)$ by 
$H_k(g)=\max_{1\le i\le k}|\<x_{*i}\,,h(g)\>|$. Let 
$$
S_{kn}(g)=\int_{G^\iy}H_k(\om_n)\,Q_g(d\om)=\int_GH_k(gg')\,\mu^n(dg'),\quad 
n=0,1,\dots\,, \tag{5.4}
$$
where $\om_n =X_n(\om)$, 
and let $S(g)=\sup_{k,n}S_{kn}(g)$. 
Note that the sequence $\{H_k\}_{k=1}^\infty$ is nondecreasing and each 
$H_k$ is $\mu$-subharmonic. Hence, $S_{kn}\le S_{k'n'}$ whenever $k\le k'$ 
and $n\le n'$. Moreover, the sequence $\{H_k\}_{k=1}^\iy$ converges 
pointwise to $H$ while for any fixed $k$, $\lim_{n\to\iy}H_k(\om_n) 
=\max_{1\le i\le k}|\<x_{*i}\,,f(\om)\>|$ 
u.a.e.. 
Hence, 
$$\gather
S(g)=\sup_k\sup_n S_{kn}(g)=\sup_k\lim_{n\to\iy}S_{kn}(g)= 
\sup_k\int_{G^\iy}\max_{1\le i\le k}|\<x_{*i}\,,f(\om)\>|\,\,Q_g(d\om)\\
=\int_{G^\iy}\sup_k\max_{1\le i\le k}|\<x_{*i}\,,f(\om)\>|\,\,Q_g(d\om)
=\int_{G^\iy}\|f(\om)\|\,Q_g(d\om),\tag{5.5}
\endgather
$$
and 
$$\gather
S(g)=\sup_n\sup_k S_{kn}(g) =\sup_n\lim_{k\to\iy} S_{kn}(g)
=\sup_n\inl_{G^\iy}H(\om_n)\,\,Q_g(d\om)\\
=\lim_{n\to\iy}\int_{G^\iy}H(\om_n)\,\,Q_g(d\om)
=\int_{G^\iy}F(\om)\,\,Q_g(d\om).\tag{5.6}
\endgather
$$ 
This completes the proof. \qed\enddemo
\remark{\bf Remark 5.2} It is a standard result of the theory of vector 
valued martingales [34, Proposition V-2-8] that when $E$ is assumed separable 
then the sequence $h\circ X_n$ converges in norm u.a.e. to $f$. This result 
can be easily deduced from the second statement of our lemma by mimicking a 
part of the argument on p.\,110 in [34]. 
\endremark

\proclaim{Proposition 5.3} $R_E$ is an equivariant isometry of \phantom{,}
$\Bbb L^\iy_i(\mu,E)$ onto\phantom{,} $\Bbb H_\mu(E)$. Moreover, for every 
$h\in\Bbb H_\mu(E)$ the sequence $h\circ X_n$ converges weakly* u.a.e. to 
$R_E^{-1}h$ and the sequence $\|h\circ X_n\|$ converges u.a.e. to 
$\|R^{-1}_Eh\|$. 
\endproclaim 

\demo{Proof}By Lemma 5.1 it suffices to prove the first statement. 

We know that $R_E$ is an equivariant contraction into $\Bbb H_\mu(E)$. Let 
$f\:G^\iy\to E$ be a universally essentially bounded invariant random 
variable and $h=R_Ef$. Using the first or the second statement of 
Lemma 5.1 it immediately follows that $\|f\|_u\le \|h\|_{\sup}$. Hence, $R_E$ 
is isometric. 
It remains to prove that $R_E$ is surjective. 

Given $h\in \Bbb H_\mu(E)$, for each $x_*\in E_*$ define a 
function $h_{x_*}\:G\to\Bbb C$ by $h_{x_*}(g)=\<x_*\,,h(g)\>$. Then 
$h_{x_*}\in\Bbb H_\mu$. Hence, by Proposition 4.1, the set 
$$
\Om_{x_*}=\bigl\{\om\in G^\iy\,;\,\{h_{x_*}(\om_n)\}_{n=0}^\iy\text{ 
converges}\,\bigr\}
$$ 
is universally conull and there exists $f_{x_*}\in \Bbb B_i$ such that for 
every $\om\in\Om_{x_*}$, $f_{x_*}(\om)=\lim_{n\to\iy}h_{x_*}(\om_n)$, 
and that $h_{x_*}=Rf_{x_*}$. A routine argument using separability of $E_*$ 
shows that 
$\Om=\bigcap_{x_*\in E_*}\Om_{x_*}$ is also universally conull. Now, for 
every $\om\in \Om$, the function $E_*\ni x_*\to f_{x_*}(\om)$ is a bounded 
linear functional on $E_*$. It follows that there 
is a $\Cal B^{(i)}$-measurable function $f\:G^\iy\to E$ such that for 
every $\om\in \Om$, $f(\om)=\w*lim_{n\to\iy}h(\om_n)$, and 
$h=R_Ef$. So $R_E$ is indeed surjective. 
\qed\enddemo
When $E$ is a $W^*$-algebra, then under the norm $\|\cdot\|_u$, $\Bbb 
L^\iy_i(\mu,E)$ is a $C^*$-algebra. As in the classical case we will 
denote by $\di$ the product 
$$
h_1\di h_2 =R_E[(R_E^{-1}h_1)(R_E^{-1}h_2)]
$$ 
in $\Bbb H_\mu(E)$. 

\proclaim{Corollary 5.4} If $E$ is a $W^*$-algebra then for every $h\in \Bbb 
H_\mu(E)$ the sequence $\{h\circ\nmb  X_n\}_{n=0}^\iy$ converges 
in the $\si$-strong* topology u.a.e. to $R_E^{-1}h$. Moreover, given 
$h_1,h_2\in \Bbb H_\mu(E)$ and $g\in G$, the sequence $[r_E(\mu^n)
(h_1h_2)](g)$ converges in the $\si$-strong* topology to $(h_1\di h_2)(g)$.
\endproclaim 
\demo{Proof} Since the predual of $E$ is separable, we may assume that $E$ 
acts in a separable Hilbert space $\frak H$. Since $\Bbb H_\mu(E)$ is a 
selfadjoint subspace of $\Bbb B(G,E)$ and each $h\in \Bbb H_\mu(E)$ is 
bounded, it suffices to prove convergence in the strong operator topology. 

Given $h\in\Bbb H_\mu(E)$ let $f$ be a representative of $R_E^{-1}h$. 
We need to show that there exists a universally conull invariant set $\Om$ 
such that for every $\om\in\Om$ and $\xi\in\frak H$, 
the sequence $h(\om_n)\xi$ converges in norm to $f(\om)\xi$. 

Note that the function $h_\xi(g)=h(g)\xi$ is an $\frak H$-valued bounded 
$\mu$-harmonic function and the function $f_\xi(\om)=f(\om)\xi$, 
a universally essentially bounded $\frak H$-valued invariant random variable 
with $h_\xi=R_{\frak H}f_\xi$. Hence, by Lemma 5.1 the set 
$$
\Om_\xi=\{\om\in G^\iy\,;\,\wlim_{n\to\infty}h(\om_n)\xi 
=f(\om)\xi\text{ and } \lim_{n\to\infty}\|h(\om_n)\xi\|=\|f(\om)\xi\|\,\}
$$ 
is universally conull. 
But by an elementary result on weak convergence in $\frak H$, 
$$
\Om_\xi=\{\om\in  G^\iy\,;\,\lim_{n\to\infty}h(\om_n)\xi=f(\om)\xi\,\}.
$$ 
Separability of $\frak H$ implies that 
$$
\Om=\bigcap_{\xi\in\frak H}\Om_\xi =\{\om\in G^\iy\,;\,f(\om)=\slim_{n\to\iy} 
h(\om_n)\}
$$ 
is also universally conull. This proves the first statement. 

To prove the second statement, given $h_1,h_2\in\Bbb H_\mu(E)$ let 
$f_1,f_2$ be representatives of $R^{-1}h_1$ and $R^{-1}h_2$, resp.. 
We know that $f_i(\om)=\slim_{n\to\iy}h_i(\om_n)$ u.a.e. and therefore using 
boundedness of the $h_i$'s we also have 
$(f_1f_2)(\om)=\slim_{n\to\iy}(h_1h_2)(\om_n)$ u.a.e.. Then a version of the 
Dominated Convergence Theorem yields, for every $g\in G$ and $\xi\in\frak H$, 
$$\gather 
(h_1\di h_2)(g)\xi=\int_{G^\iy}(f_1f_2)(\om)\xi\,\,Q_g(d\om)\\
=\lim_{n\to\infty}\int_{G^\iy}(h_1h_2)(\om_n)\xi\,\,Q_g(d\om)=
\lim_{n\to\infty}\int_{G^\iy}(h_1h_2)(gg')\xi\,\,\mu^n(dg')\\
=\vphantom{\int_G^{\iy}}\lim_{n\to\infty}[r_E(\mu^n)(h_1h_2)](g)\xi. 
\tag"{\qed}"
\endgather
$$
\enddemo
              
\example{\bf Remark 5.5}Write 
$\Cal H_\mu(E)$ for the space of equivalence classes modulo $\lam$ of 
bounded $E$-valued $\mu$-harmonic functions and $L^\iy_i(\mu,E)$ for 
$L^\iy(G^\iy,\Cal B^{(i)},Q_\lam,E)$, the space of equivalence classes 
modulo $Q_\lam$ of $Q_\lam$-essentially bounded $E$-valued invariant random 
variables. When $G$ is $\si$-compact there is no difficulty 
extending Proposition 5.3 and Corollary 5.4 to this setting, but 
we will not need it. An important difference between $\Bbb H_\mu(E)$ and 
$\Cal H_\mu(E)$, and between $\Bbb L^\iy_i(\mu,E)$ and $L^\iy_i(\mu,E)$ 
(which is 
already present in the classical case $E=\Bbb C$) is that $\Cal H_\mu(E)$ and 
$L^\iy_i(\mu,E)$ are duals of Banach spaces while $\Bbb H_\mu(E)$ and $\Bbb 
L^\iy_i(\mu,E)$ are, in general, not. Therefore, when $E$ is a $W^*$-algebra, 
$\Bbb L^\iy_i(\mu,E)$ and $\Bbb H_\mu(E)$ will be $C^*$-algebras but, in 
general, not $W^*$-algebras, while $\Cal H_\mu(E)$ and $L^\iy_i(\mu,E) = 
L^\iy_i(\mu)\oti E$ remain $W^*$-algebras.
\endexample

\head{\bf 6. Harmonic vectors}\endhead

As in Section 5, $E$ will denote a Banach space with a separable predual 
$E_*$ and $G$ a locally compact group. We will consider a representation 
$\pi$ of $G$ in $E$ which is the adjoint of a strongly continuous 
representation $\pi_*$ by isometries in $E_*$. Our goal is to obtain a 
Poisson formula for the $\mu$-harmonic vectors in $E$ (cf. Section 2) and 
relate it to the classical Poisson formula for the bounded $\mu$-harmonic 
functions. Our main result will require the assumption that $G$ be second 
countable, but the initial results of this section are valid for any locally 
compact group. 

Observe that when $x\in E$ is a $\mu$-harmonic vector, then the function 
$g\to\pi(g)x$ is a bounded $E$-valued $\mu$-harmonic function. Hence, the 
space $\Cal H_{\mu\pi}$ of the $\mu$-harmonic vectors in $E$ is 
isometrically 
isomorphic to a closed subspace of\phantom{,} $\Bbb H_\mu(E)$ and therefore 
also to a closed subspace, $\Bbb L_{\mu\pi}$, of 
$\Bbb L_i^\iy(\mu,E)$. We proceed to give a more explicit description of 
$\Bbb L_{\mu\pi}$. 

Recall that $G$ acts in a natural way on each of the function spaces $\Bbb 
B(G,E)$, $\Bbb H_\mu(E)$, $\Bbb B_i(E)$, $\Bbb L^\iy_i(\mu,E)$, 
and $L^\iy_i(\mu,E)$, and that we write $gf$ for $g\in G$ applied to a 
function $f$. Now, the representation $\pi$ induces a representation, 
$\hat\pi$, in each of the function spaces in question: $\hat\pi(g)$ 
transforms a function $f$ $\bigl($or a class of functions when $f\in\Bbb 
L^\iy_i(\mu,E)$ or $f\in L^\iy_i(\mu,E)$$\bigr)$ into the function 
$(\hat\pi(g)f)(\cdot)=\pi(g)f(\cdot)$. It is clear that $\hat\pi$ commutes 
with the natural action of $G$. 

It is easy to see that $\mu$-harmonic functions of the form $\pi(\cdot)x$, 
where $x\in\Cal H_{\mu\pi}$, are precisely those elements $h\in\Bbb H_\mu(E)$ 
which satisfy $g^{-1}h=\hat\pi(g)h$ for every $g\in G$, i.e, 
$h(gg')=\pi(g)h(g')$ for all $g,g'\in G$. Since the isomorphism $R_E$ of 
Proposition 5.3 is equivariant with respect to both the natural action of $G$ 
and the representation $\hat\pi$, it follows that $\Cal H_{\mu\pi}$ is 
isometrically isomorphic to 
$$
\Bbb L_{\mu\pi}=\{f\in \Bbb L^\iy_i(\mu,E)\,;\,g^{-1}f=\hat\pi(g)f\text{ for 
every }g\in G\,\}.
$$ 
More precisely, we have the following.

\proclaim{Theorem 6.1} The mapping $R_\pi\:\Bbb L_{\mu\pi}\to E$ given by 
$$
R_\pi f=(R_Ef)(e)=\int_{G^\iy}f(\om)\,Q_e(d\om)\tag{6.1}
$$
is an isometric isomorphism of $\,\Bbb L_{\mu\pi}$ onto $\Cal H_{\mu\pi}$. 
Furthermore, given $x\in \Cal H_{\mu\pi}$, the sequence $\pi(\om_n)x$ 
converges weakly* u.a.e. to $R_\pi^{-1}x$. 
\endproclaim 

Note that the elements of $\Bbb L_{\mu\pi}$ are equivalence classes 
modulo $\Cal N_u$ of those universally essentially bounded 
$\Cal B^{(i)}$-measurable functions $f\:G^\iy\to E$ which 
have the property that for each $g\in G$, $f(g\om)=\pi(g)f(\om)$ u.a.e.. Let 
$$
\Bbb B_{i\pi}=\{f\in \Bbb B_i(E)\,;\,f(g\om)=\pi(g)f(\om)\text{ for 
every $g\in G$ and $\om\in G^\iy$}\}. 
$$
Trivially, every element of 
$\Bbb B_{i\pi}$ represents an element of $\Bbb L_{\mu\pi}$. But it is also 
true that every element of $\Bbb L_{\mu\pi}$ can be represented by an element 
of $\Bbb B_{i\pi}$. 
\proclaim{Corollary 6.2}$\Bbb L_{\mu\pi}$ is the quotient of \phantom,$\Bbb 
B_{i\pi}$ modulo $\Cal N_u$. 
\endproclaim 
\demo{Proof}It suffices to show that if $f\:G^\iy\to E$ is a universally 
essentially bounded invariant random variable such that for each $g\in G$, 
$f(g\om)=\pi(g)f(\om)$ u.a.e., then there exists $f'\in \Bbb B_{i\pi}$ with 
$f=f'$ u.a.e.. But with $x=R_\pi f$, the set $\{\om\nmb\in\nmb 
G^\iy\nmb\,\nmb ;\,\w*lim_{n\to\iy}\pi(\om_n)x\text{ exists}\}$ is obviously 
$G$-invariant and by Proposition 6.1 it is universally conull. Define 
$f'\:G^\iy\to E$ by 
$$
f'(\om)=\cases \w*lim_{n\to\iy}\pi(\om_n)x,&\text{when the limit exists}\\
0,&\text{otherwise}.\endcases\tag"\qed"
$$\enddemo
Eq.\,(6.1) can be viewed as a Poisson formula for the $\mu$-harmonic vectors. 
However, the use of the subspace $\Bbb L_{\mu\pi}$ and the norm $\|\cdot\|_u$ 
has certain serious disadvantages. 

Firstly, $\Cal H_{\mu\pi}$ is the dual of 
the quotient $E_*/J_{\mu\pi}$, because $\Cal H_{\mu\pi} 
=J_{\mu\pi}^{\,\perp}$. So $\Bbb L_{\mu\pi}$ is also the dual of 
$E_*/J_{\mu\pi}$. But, in general, $\Bbb L_i^\iy(\mu,E)$ has no predual and 
therefore it is not 
clear how the predual $E_*/J_{\mu\pi}$ of $\Bbb L_{\mu\pi}$ is related to the 
Poisson formula and the boundary space $(G^\iy,\Cal B^{(i)})$. 

Secondly, as 
we already pointed out in Section 4, when $G$ is second countable then 
for the purpose of representing the classical $\mu$-harmonic functions, the 
badly behaved Borel space $(G^\iy,\Cal B^{(i)})$ can be replaced by a 
standard Borel $G$-space or even a topological $G$-space. 
Recall that this regularization of the boundary relies on the 
fact that $L^\iy_i(\mu)$ is an abelian $W^*$-algebra. 

Now, 
let $\be$ be any fixed finite measure equivalent to $\lam$. Then 
$L^\iy_i(\mu,E)$ is canonically the dual of $L^1_i(\mu,E_*)=L^1(G^\iy,\Cal 
B^{(i)},Q_\be,E_*)$, the space of equivalence classes of $Q_\be$-integrable 
$\Cal B^{(i)}$-measurable functions $f\:G^\iy\to E_*$. 
As we will see, the 
Poisson formula for the $\mu$-harmonic vectors becomes very satisfactory once 
the subspace $\Bbb L_{\mu\pi}\sub \Bbb L^\iy_i(\mu, E)$ is replaced by the 
weakly* closed subspace 
$$
L_{\mu\pi}=\{f\in L^\iy_i(\mu,E)\,;\,g^{-1}f 
=\hat\pi(g)f\,\text{ for every }g\in G\,\}\sub L^\iy_i(\mu,E).
$$ 

{\bf For the remainder of this article we will assume that $G$ is 
second countable.} 

In the setting of Example 2.7 one readily obtains:

\proclaim{Proposition 6.3} $L_{\mu\pi}$ equals the crossed product 
$L^\iy_i(\mu)\ti_{\pi_l}G$ where $\pi_l$ is the left regular 
representation in $L^\iy_i(\mu)$. 
\endproclaim 

\proclaim{Lemma 6.4} Let $f\:G^\iy\to E$ be a bounded invariant random 
variable such that for each $g\in G$, $f(g\om)=\pi(g)f(\om)$ for 
$Q_\lam$-a.e. $\om\in G^\iy$. Then there exists an 
invariant random variable $f'\:G^\iy\to E$ such that $f=f'$ 
$\,Q_\lam$-a.e.,
$\|f'\|_{\sup}\le\|f\|_{\sup}$, and $\pi(g)f'(\om)=f'(g\om)$ for every $g\in 
G$ 
and $\om\in G^\iy$. 
\endproclaim 
\demo{Proof}Consider the Borel $G$-space $(G^\iy,\Cal B^\iy)$. As a countable 
product of standard Borel spaces, $(G^\iy,\Cal B^\iy)$ is standard. Moreover, 
$Q_\lam$ is a $\si$-finite invariant measure on $(G^\iy,\Cal B^\iy)$ and the 
mapping $G\ti G^\iy\ni (g,\om)\to\nmb g\om\in G^\iy$ is Borel. 

Next, let 
$E_f$ denote the closed ball $\bar B_E(0,\|f\|_{\sup})$ of radius 
$\|f\|_{\sup}$ and centre $0$ in $E$. With the $G$-action given by $\pi$ and 
with the weak* topology, $E_f$ is a topological $G$-space and the mapping 
$G\times E_f\ni (g,x)\to \pi(g)x\in E_f$ is continuous. 
Therefore this mapping is also Borel (with respect to the product Borel 
structure on $G\ti E_f$). Thus $E_f$ is a standard Borel $G$-space with the 
mapping $G\times E_f\ni (g,x)\to \pi(g)x\in E_f$ Borel. $f$ is a Borel 
function of $G^\iy$ into $E_f$ such that for every $g\in G$, 
$f(g\om)=\pi(g)f(\om)$ $Q_\lam$-a.e.. We are therefore in a position to apply 
[45, Proposition B.5, p.\,198] to conclude that there exists 
a $\Cal B^\iy$-measurable function $f''\:G^\iy\to E_f$ 
such that $f''=f$ $Q_\lam$-a.e. and $f''(g\om)=\pi(g)f''(\om)$ for every 
$g\in G$ and $\om\in G^\iy$. 

Let 
$$\gather
\G= \{\om\in G^\iy\,;\,f''(\om)=f''(\vt^n(\om))\text{ for every }n\ge 
0\}\\=\bigcap_{n=0}^\iy\vt^{-n}\bigl(\{\om\in\nmb G^\iy\,;\,\alb f''(\om) 
=f''(\vt(\om))\}\bigr) \endgather 
$$ 
where $\vt$ is the Markov shift, cf. Section 4. 
Note that $\G\sub\vt^{-1}(\G)$, $\G$ is $G$-invariant, $\G\in\Cal B^\iy$, and 
$Q_\lam(G^\iy-\G)=0$ (the latter because $\vt$ preserves the measure class of 
$Q_\lam$, cf. Eq.\,(4.5)). 

Let 
$$
\Del=\bigcap_{n=0}^\iy\vt^{-n}(G^\iy-\G)=\{\om\in 
G^\iy-\G\,;\,\vt^n(\om)\in G^\iy-\G\text{ for every } n\ge 1\}.
$$ 
Then $\vt^{-1}(\Del)=\Del$, i.e., $\Del\in\Cal B^{(i)}$, and 
$\Del$ is $G$-invariant. Note also that if $\om\in (G^\iy-\G)\cap 
(G^\iy-\Del)$ then there exists $k(\om)\ge 1$ such that 
$\vt^{k(\om)}(\om)\in\G$, and if $\vt^i(\om)\in\G$ and 
$\vt^{j}(\om)\in\G$ for some $i,j\ge 1$, then 
$f''(\vt^i(\om))=f''(\vt^j(\om))$. 

Define a function $f'\: G^\iy\to E$ by
$$f'(\om)=\cases 0,&\text{for $\om\in\Del$}\\
         f''(\om),&\text{for $\om\in\G$}\\
         f''(\vt^{k(\om)}\om),&\text{for } \om\in G^\iy-(\G\cup \Del). 
\endcases 
$$
It is straightforward to verify that $f'$ is an invariant random variable 
equal to $f$ $Q_\lam$-a.e. and satisfying $\|f'\|_{\sup}\le \|f\|_{\sup}$ 
as well as $f'(g\om)=\pi(g)f'(\om)$ for every $g\in G$ and $\om\in G^\iy$. 
\qed\enddemo
       
\proclaim{Lemma 6.5}Given $f_1,f_2\in \Bbb B_{i\pi}$ the following conditions 
are equivalent\rom :
\roster
\item"{(i)}" $f_1=f_2$ u.a.e.,
\item"{(ii)}" $f_1=f_2$ $Q_\lam$-a.e.,
\item"(iii)" $f_1=f_2$ $Q_e$-a.e..
\endroster
Moreover, for every $f\in \Bbb B_{i\pi}$, $\|f\|$ is constant u.a.e.. 
\endproclaim 
\demo{Proof}Note that the set $\{\om\in G^\iy\,;\,f_1(\om)=f_2(\om)\}$ is 
$G$-invariant. The equivalence follows immediately from this observation 
and Eqs\,(4.1) and (4.2). The second statement follows from Corollary 4.4 and 
the observation that the function $F(\om)=\|f(\om)\|$ is $G$-invariant. 
\qed\enddemo

By Corollary 6.2 and Lemmas 6.4 and 6.5 there exists a surjective 
mapping $\Phi\:\Bbb L_{\mu\pi}\to L_{\mu\pi}$ such that $\Phi [f]_u 
=[f]_{Q_\lam}$ for every $f\in \Bbb B_{i\pi}$ where $[\cdot]_u$ and 
$[\cdot]_{Q_\lam}$ denote equivalence classes modulo $\Cal N_u$ and $Q_\lam$, 
resp.. It is clear that $\Phi$ is a linear mapping. It is also an 
isometry because by Lemma 6.5 the essential sup-norm of $f\in\Bbb B_{i\pi}$ 
modulo $\Cal N_u$ is the same as the essential sup-norm modulo 
$Q_\lam$-null sets. It follows that there exists an isometry 
$R_\pi$ of $L_{\mu\pi}$ onto $\Cal H_{\mu\pi}$ such that 
$$
R_\pi [f]_{Q_\lam}= (R_Ef)(e)=\int_{G^\iy}f(\om)\,Q_e(d\om)\tag{6.2}
$$
for every $f\in\Bbb B_{i\pi}$. 

\proclaim{Theorem 6.5} $R_\pi$ is an isometric weak* homeomorphism of 
$L_{\mu\pi}$ onto $\Cal H_{\mu\pi}$. 
\endproclaim 

\demo{Proof}It remains to prove that $R_\pi$ and $R_\pi^{-1}$ are weakly* 
continuous. By Lemma 4.8 it is enough to prove this for $R_\pi$. 
Using the Krein-Smulian theorem, it suffices to show that the restriction of 
$R_\pi$ to the closed unit ball of $L_{\mu\pi}$ is weakly* continuous at 
zero. Let $F_j$ be a net in the closed unit ball of $L_{\mu\pi}$, converging 
weakly* to $0$. It is enough to show that $\lim_{j}\<x_*\,,R_\pi F_j\>=0$ 
for all $x_*$ in a dense subset $D$ of $E_*$. 

Using an approximate identity in $L^1(G)\sub M(G)$ one can see that the set 
$$
D=\bigl\{\int_G\pi_*(g)x_*\,\,\si (dg)\,;\,\si\in L^1(G),\ x_*\in E_*\bigr\}
$$ 
is dense in $E_*$. Using Lemma 6.4, for each $i$ choose a 
representative $f_j$ of $F_j$ in $\Bbb B_{i\pi}$. Let 
$y_*=\int_G\pi_*(g)x_*\,\si (dg)\in D$. Then 
$$\gather 
\<y_*\,,R_\pi F_j\>=\int_{G^\iy}\<y_*\,,f_j(\om)\>\,Q_e(d\om)
=\int_{G^\iy}\biggl(\int_G \<\pi_*(g)x_*\,,f_j(\om)\>\,\si(dg) 
\biggr)\,Q_e(d\om)\\
=
\int_{G^\iy}\biggl(\int_G \<x_*\,,\pi(g^{-1})f_j(\om)\>\,\si(dg) 
\biggr)\,Q_e(d\om)\\
=\int_{G^\iy}\biggl(\int_G \<x_*\,,f_j(g^{-1}\om)\>\,\si(dg) 
\biggr)\,Q_e(d\om)
=
\int_{G^\iy}\<x_*\,,f_j(\om)\>\,Q_{\tilde\si}(d\om). \tag{6.3}
\endgather
$$
Let $\be$ be a finite measure equivalent to $\lam$. Then $Q_{\tilde\si}\ll 
Q_\be$. Let $s$ be a version of the Radon-Nikodym derivative of the 
restriction of $Q_{\tilde\si}$ to $\Cal B^{(i)}$ with respect to the 
restriction of $Q_\be$ to $\Cal B^{(i)}$. Then 
$$
\<y_*\,,R_\pi F_j\>=
\int_{G^\iy}\<s(\om)x_*\,,f_j(\om)\>\,\,Q_\be(d\om).
$$
Since the function $\om\to s(\om)x_*$ defines an element of $L^1_i(\mu,E_*) 
=L^\iy_i(\mu,E)_*$, and $\w*lim_jF_j=0$, it becomes clear that 
$\lim_j\<y_*\,,R_\pi F_j\>=0$. 
\qed\enddemo

\remark{\bf Remark 6.7} Let $(\Cal X,\Cal A,\al)$ be a $\mu$-boundary which 
is a standard Borel space with the mapping $G\ti\Cal X\ni (g,t)\to gt\in\Cal 
X$ Borel, and let $\Phi$ be the equivariant isomorphism of $L^\iy(\Cal X,\Cal 
A,\al)$ onto $L^\iy_i(\mu)$. 
It can be shown [19] that there exists a $\Cal B^{(i)}$-measurable function 
$F\:G^\iy\to\Cal X$ which induces $\Phi$ in the sense that $\Phi f=f\circ F$, 
$f\in L^\iy(\Cal X,\Cal A,\al)$, and is itself equivariant when restricted 
to a suitable $G$-invariant subset $\Om\in\Cal B^{(i)}$. The measure 
$\rho=FQ_e$ is the Poisson kernel on the $\mu$-boundary $(\Cal X,\Cal 
A,\al)$: 
$$
(R\Phi f)(g)=\int_{\Cal X}f(gt)\,\rho(dt)\pmod\lam
$$ 
for every $f\in 
L^\iy(\Cal X,\Cal A,\al)\,$.  
\footnote{\ This is shown in [19] for so-called continuous $\mu$-boundaries. 
The result for standard $\mu$-boundaries requires a little extra work 
relying on the theory of pointwise realizations of $L^\iy$-spaces 
and homomorphisms between them [28].} 

Using the mapping $F$ one can obtain a version of the Poisson formula 
(6.2) on the standard $\mu$-boundary $(\Cal X,\Cal A,\al)$. This can be 
particularly useful when a ``natural'' realization of the $\mu$-boundary is 
known. For example, when $G$ is almost connected and $\mu$ is spread out, 
then the natural pointwise realization of the $\mu$-boundary turns out to be 
a homogeneous space $G/H$ of $G$ [23]; $\Cal H_{\mu\pi}$ is then isomorphic 
to the space of equivariant maps from $G/H$ to $E$, modulo the family of null 
sets of $G/H$.
\endremark
              
\remark{\bf Remark 6.8} Theorem 6.6 applies, in particular, to the classical 
Choquet-Deny equation in $E=L^\iy(G)$. In this case $L_{\mu\pi}$ must 
therefore be isomorphic to $L^\iy_i(\mu)$. 
The isomorphism $\Psi= R_\pi^{-1}R\: L^\iy_i(\mu)\to L_{\mu\pi}$ 
(where $R$ is as in Proposition 4.5) can be described explicitly as follows. 
Given $f\in \Bbb B_i$, let $\dot f$ denote the function $\dot f\:G^\iy\times 
G\to\Bbb C$ given by $\dot f(\om,g)=f(g\om)$. Using the fact that the mapping 
$G^\iy\times G\ni (\om,g)\to g\om\in G^\iy$ is measurable with respect to the 
$\si$-algebras $\Cal B^\iy\times \Cal B$ and $\Cal B^\iy$, it follows that 
for a fixed $\om$, $f(\om,\cdot)\in \Bbb B(G)$, and for every $\nu\in 
L^1(G)\sub M(G)$, the function $\om\to \int_G\dot f(\om,g)\,\nu(dg)\in\Bbb C$ 
is $\Cal B^{(i)}$-measurable. Hence, if we define $\ddot f\:G^\iy\to 
L^\iy(G)$ by $\ddot f(\om) =[\dot f(\om,\cdot)]_\lam$, then $\ddot f\in \Bbb 
B_i\bigl(L^\iy(G)\bigr)$. It is immediate that $\ddot f\in \Bbb 
B_{i\pi}\bigl(L^\iy(G)\bigr)$, so that $[\ddot f]_{Q_\lam}\in L_{\mu\pi}$. 
Moreover, $R_\pi [\ddot f]_{Q_\lam}=R[f]_{Q_\lam}$. Therefore $\Psi$ is given 
by $\Psi [f]_{Q_\lam}=[\ddot f]_{Q_\lam}$ for every $f\in \Bbb B_i$. 
We also note that $\Psi$ is equivariant with respect to the natural actions 
of $G$ on $L^\iy_i(\mu)$ and $L_{\mu\pi}$ (the left regular representations). 
\endremark

\remark{\bf Remark 6.9} Let $\Cal F_{\mu\pi}\sub \Cal H_{\mu\pi}$ denote the 
subspace of trivial solutions of the Choquet-Deny equation in $E$. 
Recall that $x\in \Cal F_{\mu\pi}$ if and only if $\pi(g)x=x$ for every $g\in 
G_\mu$. Let $F_{\mu\pi}=R_\pi^{-1}\Cal F_{\mu\pi}$. We claim that given $f\in 
\Bbb B_{i\pi}$, $[f]_{Q_\lam}\in F_{\mu\pi}$ if and only if $f$ is constant 
$Q_e$-a.e.. 

Indeed, necessity follows from the last statement of Theorem 6.1: 
if $x=R_\pi [f]_{Q_\lam}$ then 
$$
f(\om)=\w*lim_{n\to\iy}\pi(\om_n)x\quad \text{$Q_e$-a.e..} 
$$
But 
$$
\int_{G^\iy}\|\pi(\om_n)x-x\|\,Q_e(d\om)=\int_G\|\pi(g)x-x\|\,\mu^n(dg)=0.
$$ 
So for each $n$, $\pi(\om_n)x=x$ $Q_e$-a.e. and therefore $f(\om)=x$ 
$Q_e$-a.e.. Conversely, suppose that there is a set $\Om\in \Cal B^{(i)}$ 
with $Q_e(\Om)=1$ and $f(\om)=x$ for every $\om\in\Om$. Using Eqs\,(4.7) and 
(4.2) we obtain that $Q_g(\Om)=1$ for $\mu$-a.e. $g\in G$. Hence, 
$$
\pi(g)x=\int_{G^\iy}\pi(g)f(\om)\,Q_e(d\om)= \int_{G^\iy}f(g\om)\,Q_e(d\om) 
=\int_{G^\iy}f(\om)\,Q_g(d\om)=x
$$ 
for $\mu$-a.e. $g\in G$. Since the mapping 
$g\to\pi(g)x$ is continuous with respect to the weak* topology on $E$, this 
implies that $\pi(g)x=x$ for every $g\in G_\mu$. 

The above description of $F_{\mu\pi}$ simplifies when 
$G=G_\mu$ because in this case the constant function $f(\om)=x$, $x\in 
\Cal F_{\mu\pi}$, belongs to $\Bbb B_{i\pi}$. Therefore in this case 
$F_{\mu\pi}$ consists of equivalence classes modulo $Q_\lam$ of such 
constant functions. 

For a general $\mu$, $\Cal H_{\mu\pi}=\Cal F_{\mu\pi}$ if and only if every 
$f\in \Bbb B_{i\pi}$ is constant modulo $Q_e$. This will be the case whenever 
the classical Choquet-Deny equation in $L^\iy(G)$ has only trivial solutions. 
However, recall that for a certain type of representations we have $\Cal 
H_{\mu\pi}=\Cal F_{\mu\pi}$, regardless of what $\Cal H_\mu$ is (cf. 
Proposition 2.2, Corollaries 2.3 and 2.4, and Proposition 2.5). Thus for 
certain type of representations it is always true that every bounded $\Cal 
B^{(i)}$-measurable equivariant function $f\:G^\iy\to E$ is constant modulo 
$Q_e$ (and modulo $Q_\lam$ when $G=G_\mu$). This seems to be an interesting 
and completely unexplored property of the $\mu$-boundaries. 
\endremark

\example{Example 2.6 (continued)} The solutions of the Choquet-Deny equation 
in $E=M(\Cal X)$ are in one-to-one correspondence with equivalence 
classes modulo $Q_\lam$ of bounded $\Cal B^{(i)}$-measurable equivariant 
functions $f\:G^\iy\to M(\Cal X)$. Suppose that $\si =R_\pi \vp$ is a 
$\mu$-stationary 
measure, i.e., $\si\in \Cal H_{\mu\pi}\cap M_1(\Cal X)$. Using the last 
statement of Theorem 6.1 it is easy to see that $\vp(\om)\in M_1(\Cal 
X)\phantom i$ $Q_\lam$-a.e.. More precisely, $\vp=[f]_{Q_\lam}$ where 
$f\in \Bbb B_{i\pi}$ is such that $f(\om)\in 
M_1(\Cal X)$ for every $\om$ in a universally conull $G$-invariant 
set $\Om \in \Cal B^{(i)}$. Let $\Bbb B_{i\pi 1}$ denote the subset of $\Bbb 
B_{i\pi}$ consisting of such functions. Thus $\mu$-stationary measures are in 
one-to-one correspondence with equivalence classes of the elements of 
$\Bbb B_{i\pi1}$. 

In general, the set of $\mu$-stationary measures, and 
therefore also $\Bbb B_{i\pi 1}$, can be empty (cf. Corollary 2.3). It 
is never empty when $\Cal X$ is compact. In this case the 
existence of an equivariant function $f\in \Bbb B_{i\pi 1}$ can be also 
deduced using [45, Proposition 4.3.9] and the fact that the action of $G$ on 
every $\mu$-boundary is amenable. The fact that the existence of a 
$\mu$-stationary measure on a locally compact second countable $G$-space 
implies the existence of an equivariant function $f\in\Bbb B_{i\pi1}$ was 
first observed and exploited by Furstenberg [10]. \ee 
\endexample 

Suppose that $E$ is a $W^*$-algebra and $\pi(g)\in \Aut(E)$ for every $g\in 
G$. Then $L_{\mu\pi}$ is a weakly* closed *-subalgebra of 
$L^\iy_i(\mu,E) = L^\iy_i(\mu)\oti E$, and so both $L_{\mu\pi}$ and $\Cal 
H_{\mu\pi}$ are themselves $W^*$-algebras when the multiplication, 
$\di$, in $\Cal H_{\mu\pi}$ is given by 
$$
x_1\di x_2 
=R_\pi[(R_\pi^{-1}x_1)(R_\pi^{-1}x_1)].
$$ 
A simple computation and Corollary 5.4 yield: 

\proclaim{Corollary 6.10}If $E$ is a $W^*$-algebra and $\pi(g)\in \Aut(E)$ 
for every $g\in G$, then for every $x\in \Cal H_{\mu\pi}$ the sequence 
$\pi(\om_n)x$ converges $Q_\lam$-a.e. in the $\si$-strong* topology  
to $R_\pi^{-1}x$. Moreover, given $x_1,x_2\in \Cal H_{\mu\pi}$, the sequence 
$\pi(\mu^n)(x_1x_2)$ converges in the $\si$-strong* topology to $x_1\di x_2$. 
\endproclaim
\flushpar
Evidently, the product $\di$ coincides with the product of $E$ whenever the 
Choquet-Deny theorem is true for $\mu$. We note that in the setting of 
Example 2.7, since $\Cal H_{\mu\pi}$ contains a copy of $\Cal H_\mu$, the 
product $\di$ coincides with the product of $B(L^2(G))$ if and only if the 
Choquet-Deny theorem is true. 
            
\head{\bf 7. Approximation properties in the predual and projections}\endhead

The goal of this concluding section is to present some additional 
properties of the space of $\mu$-harmonic vectors and the Poisson 
formula. We will continue to work in the setting of Section 6. Recall that 
$G$ is assumed second countable. 

We begin with a few results which link the predual $L_{\mu\pi*}$ of 
$L_{\mu\pi}$ to the predual of $E$. Let $R_{\pi*}\:E_*\to L_{\mu\pi*}$ 
denote the preadjoint of $R_\pi$ and $(E_*)_1$ and $(L_{\mu\pi *})_1$ 
the unit balls in $E_*$ and $L_{\mu\pi*}$, resp.. Our first result summarizes 
general properties of weakly* continuous isometries between dual Banach 
spaces. 

\proclaim{Theorem 7.1} $R_{\pi *}E_*= L_{\mu\pi*}$ and 
$R_{\pi *}\bigl((E_*)_1\bigr)$ is norm dense in $(L_{\mu\pi *})_1$. 
\endproclaim

\proclaim{Theorem 7.2} If $E$ is a W*-algebra and $\pi(g)\in\Aut(E)$ for 
every $g\in G$, then $R_{\pi *}N_E$ is norm dense in $N_{L_{\mu\pi}}$, where 
$N_E$ and $N_{L_{\mu\pi}}$ denote the sets of normal states on $E$ and 
$L_{\mu\pi}$, resp..
\endproclaim 
\demo{Proof}Since $R_\pi$ preserves positivity and maps the unit of 
$L_{\mu\pi}$ to the unit of $E$, it is clear that $R_{\pi *}N_E\sub 
N_{L_{\mu\pi}}$. The density follows by a routine application of the 
Hahn-Banach theorem. 
\qed\enddemo

Recall that in the case of the classical Choquet-Deny equation in $L^\iy(G)$, 
$R_{\pi*}\vp= \vp*Q_e$ for every $\vp\in L^1(G)\sub M(G)$. So Theorem 7.2 
says that probability measures of 
the form $\vp*Q_e$, where $\vp$ is a probability measure in $L^1(G)\sub 
M(G)$, are norm dense in the set of probability measures in $L^1_i(\mu)$. 
When $\mu$ is spread out, it can be shown that this is equivalent to the 
condition that the convex hull of the orbit of $Q_e$ under the natural 
action of $G$ be dense in the set of probability measures in $L^1_i(\mu)$. 
An abstraction of this property of the action of $G$ on the $\mu$-boundary 
is the concept of a strongly approximately transitive action which plays an
important role in the theory of $\mu$-boundaries [20-23]. When $\mu$ is 
not necessarily spread out, Theorem 7.2 implies that the action of $G$ on the 
$\mu$-boundary is approximately transitive, a concept first introduced by 
Connes and Woods in the context of the theory of Neumann algebras [6]. 

Our next result is a generalization of Th\'eor\`eme 1 of Derriennic [9]. 

\proclaim{Theorem 7.3}If $x_*\in E_*$ then 
$$
\|R_{\pi*}x_*\|=\inf_{n\ge 
1}\bigl\|\tfrac1n\sum_{i=1}^n\pi_*(\tilde\mu^i)x_*\bigr\| =\lim_{n\to\iy}
\bigl\|\tfrac1n\sum_{i=1}^n\pi_*(\tilde\mu^i)x_*\bigr\|. 
$$
\endproclaim 
\demo{Proof}The sequence $\|\sum_{i=1}^n\pi_*(\tilde\mu^i)x_*\bigr\|$, 
$n=1,2,\dots$ is subadditive. Hence, the second equality is elementary. 
Moreover, from the identity $\pi(\mu)R_\pi=R_\pi$ it immediately follows that 
$$
\|R_{\pi*}x_*\|\le\lim_{n\to\iy}\bigl\|\tfrac1n\sum_{i=1}^n\pi_*(\tilde\mu
^i)x_*\bigr\|.
$$ 
It remains to establish the opposite inequality. 

Let $\vep > 0$ be given. Then for every $n\ge 1$ we can find $x_n\in E$ with 
$\|x_n\|\le 1$ and 
$$
|\<x_*\,,\tfrac1n\sum_{i=1}^n\pi(\mu^i)x_n\>|
=|\<\tfrac1n\sum_{i=1}^n\pi_*(\tilde\mu^i)x_*\,,x_n\>|\ge 
\bigl\|\tfrac1n\sum_{i=1}^n\pi_*(\tilde\mu^i)x_*\bigr\|-\vep.
$$
Define $z_n$ by $z_n=\frac1n\sum_{i=1}^n\pi(\mu^i)x_n$. Clearly, $z_n$ is a 
sequence in the closed unit ball of $E$. Since this ball is weakly* compact 
there is a convergent subnet $z_{n_j}$ and, obviously, $z=\w*lim_jz_{n_j}\in 
\Cal H_{\mu\pi}$. Hence, there is $\zeta$ in the closed unit ball of 
$L_{\mu\pi}$ with $z=R_\pi\zeta$. Thus
$$\gather 
\lim_{n\to\iy}
\bigl\|\tfrac1n\sum_{i=1}^n\pi_*(\tilde\mu^i)x_*\bigr\|-\vep=\lim_j
\bigl\|\tfrac1{n_j}\sum_{i=1}^{n_j}\pi_*(\tilde\mu^i)x_*\bigr\|-\vep\\ \le 
\lim_j\bigl|\bigl\<x_*\,,\tfrac1{n_j}\sum_{i=1}^{n_j}\pi(\mu^i)x_{n_j}
\bigr\>\bigr|=|\<x_*\,,R_\pi\zeta\>|=|\<R_{\pi *}x_*\,,\zeta\>| \le\|R_{\pi 
*}x_*\|.
\endgather
$$
Since $\vep$ is arbitrary, we are done. \qed\enddemo

We will now prove the existence of a norm 1 projection of $E$ onto 
$\Cal H_{\mu\pi}$ and explore some of the consequences of this result. 

\proclaim{Theorem 7.4}If\phantom{,} $\Cal H_{\mu\pi}\ne \{0\}$ then there 
exists a projection, $K$, of norm 1, of $E$ onto 
$\Cal H_{\mu\pi}$, which commutes with every weakly* continuous linear 
operator $T\:E\to E$ commuting with $\pi$. If $E$ is a $W^*$-algebra and 
$\pi(g)\in\Aut(E)$ for every $g\in G$, then $K$ can be chosen 
completely positive. 
\endproclaim 
\demo{Proof}Let $\{n_j\}$ be an ultranet in $\Bbb N$ with $n_j\to\iy$. 
Since each closed ball in $E$ is weakly* compact, the ultranet 
$\frac1{n_j}\sum_{i=1}^{n_j}\pi(\mu^i)x$ converges in the weak* topology, for 
every $x\in E$. We can thus define $K\:E\to E$ by $Kx=\w*lim_j 
\frac1{n_j}\sum_{i=1}^{n_j}\pi(\mu^i)x$. It is easy to see that $K$ 
has all the desired properties. To prove complete positivity use [35, Theorem 
6.5]. \qed\enddemo

\example{\bf Remark 7.5} Let $K$ denote the projection constructed in the 
proof of Theorem 7.4. When $G_\mu$ is compact, by a version of the Ito-Kawada 
theorem [39, Theorem 2, p.\,138], the sequence $\frac1n\sum_{i=1}^n\mu^i$ 
converges weakly* to the normalized Haar measure $\om_{G_\mu}$ of $G_\mu$; 
therefore when $G_\mu$ is compact then $K=\pi(\om_{G_\mu})$. On the other 
hand, when $G_\mu$ is not compact then Lemma 2.1 implies that the kernel of 
$K$ will always contain the norm closed $\pi$-invariant subspace $E_0$ 
consisting of those vectors $x\in E$ for which the function $g\to 
\<x_*\,,\pi(g)x\>$ belongs to $C_0(G)$ for every $x_*\in E_*$. 
\endexample
                 
\example{\bf Example 7.6} Consider the classical Choquet-Deny equation in 
$E=L^\iy(G)$. It is clear from Remark 7.5 that the projection $K$ constructed 
in the proof of Theorem 7.4 is weakly* continuous if and only if $G_\mu$ is 
compact, and that when $G_\mu$ is not compact then 
$K\bigl(C_0(G)\bigr)=\{0\}$. In fact, when $G_\mu$ is not compact then not 
only the projection constructed in the proof of Theorem 7.4, but any bounded 
projection of $L^\iy(G)$ onto $\Cal H_\mu$ which commutes with every weakly* 
continuous linear operator $T\:L^\iy(G)\to L^\iy(G)$ commuting with $\pi$, 
must necessarily vanish on $C_0(G)$ and therefore cannot be weakly* 
continuous. 

Let us prove a slightly stronger result: Let $LUC(G)$ denote the set of 
bounded left uniformly continuous functions $f\:G\to\Bbb C$ and 
let $S\:LUC(G)\to \Cal H_\mu$ be 
any bounded linear operator commuting with the left regular representation. 
Then $S\bigl(LUC(G)\bigr)\sub LUC(G)$, in particular, 
$S\bigl(C_0(G)\bigr)\sub LUC(G)$. 
Hence, the mapping $C_0(G)\ni f\to (Sf)(e)$ is a well defined 
bounded linear functional on $C_0(G)$, and so there exists $\si\in M(G)$ such 
that $(Sf)(e)=\int_Gf\,d\si$ for every $f\in C_0(G)$. Using the assumption 
that $S$ commute with the left regular representation, it follows that $S\r\, 
C_0(G)=\pi(\si)\r\,C_0(G)$ (where $\pi$ is the right regular representation). 
Next, since for every $f\in C_0(G)$, 
$\pi(\mu)Sf=Sf$, i.e., $\pi(\mu*\si)f=\pi(\si)f$, we obtain $\mu*\si=\si$, 
which is a version of the Choquet-Deny equation in $M(G)$ where the left 
instead of right regular representation is used. When $G_\mu$ is not compact, 
Corollary 2.3 yields $\si=0$ and so $S$ vanishes on $C_0(G)$. Summarizing: 

\proclaim{Proposition 7.7} If $S\:LUC(G)\to\Cal H_\mu$ is a bounded linear 
operator commuting with the left regular representation then $S$ vanishes on 
$C_0(G)$. \ee 
\endproclaim 
\endexample

\example{\bf Example 2.7 (continued)} Let us identify the elements of 
$L^\iy(G)$ with the corresponding multiplication operators in $L^2(G)$. 
Moreover, let 
$$\align 
B_0(L^2(G))=\{&A\in B(L^2(G))\,;\, \text{for every }T\in\Cal T(L^2(G))\\
&\text{ the function }g\to \tr[T(\pi(g)A)]\text{ vanishes at}\mb\text{ 
infinity}\}. 
\endalign
$$ 
$B_0(L^2(G))$ is a norm closed subspace containing the ideal of 
compact operators, as well as $C_0(G)$. 

\proclaim{Proposition 7.8} When $G_\mu$ is not compact then every 
bounded linear operator $S\: B(L^2(G))\to\Cal H_{\mu\pi}$ which commutes 
with every weakly* continuous operator commuting with $\pi$, vanishes on 
$B_0(L^2(G))$. In particular, any bounded projection of $B(L^2(G))$ onto 
$\Cal H_{\mu\pi}$ which commutes with every weakly* continuous operator 
commuting with $\pi$, vanishes on $B_0(L^2(G))$ and therefore cannot be 
weakly* continuous (in fact, is even singular).
\endproclaim 
\demo{Proof} 
Given $\xi,\eta\in L^2(G)$, let $\Psi_{\xi\eta}\: B(L^2(G))\to LUC(G)\sub 
B(L^2(G))$ denote the linear mapping which associates with each $A\in 
B(L^2(G))$ the multiplication operator defined by the left uniformly 
continuous function $h_{\xi\eta A}(g)=\<(\pi(g)A)\xi\,,\eta\>$. It is 
easy to see that when $A_\al$ is a norm bounded net in $B(L^2(G))$ which 
converges to zero weakly*, then the net $h_{\xi\eta A_\al}$ converges to zero 
uniformly on compacta. This implies, via a routine argument, that 
$\Psi_{\xi\eta}$ is weakly* continuous. It is also straighforward to see that 
$\Psi_{\xi\eta}$ commutes with $\pi$. 

Let $e_j$ be a net of positive integrable functions on $G$ forming a bounded 
approximate identity in $L^1(G)$. 
Since for a given $f\in LUC(G)$, $\Psi_{\sqrt{e_j}\sqrt{e_j}}f=e_j^\star *f$, 
where $^\star$ denotes the involution in $L^1(G)$, it follows that the net 
$\Psi_{\sqrt{e_j}\sqrt{e_j}}\r LUC(G)$ converges in the strong operator 
topology to the identity operator on $LUC(G)$. Now, if $T\in 
B\bigl(B(L^2(G))\bigr)$ commutes with every weakly* continuous operator 
commuting with $\pi$ then $T\Psi_{\sqrt{e_j}\sqrt{e_j}} 
=\Psi_{\sqrt{e_j}\sqrt{e_j}}T$, and it follows that $T\bigl(LUC(G)\bigr)\sub 
LUC(G)$. 
This applies, in particular, to our operator $S\: B(L^2(G))\to\Cal 
H_{\mu\pi}$. Therefore 
$$
S\bigl(LUC(G)\bigr)\sub LUC(G)\cap \Cal H_{\mu\pi}\sub \Cal 
H_\mu.
$$ 
Thus using Proposition 7.7 we conclude that $S\bigl(C_0(G)\bigr)=\{0\}$. But 
$\Psi_{\xi\eta}\bigl(B_0(L^2(G))\bigr)\mb\sub C_0(G)$ for 
each $\xi,\eta\in L^2(G)$. Hence, 
$$
\Psi_{\xi\eta}S\bigl(B_0(L^2(G))\bigr) 
=S\Psi_{\xi\eta}\bigl(B_0(L^2(G))\bigr)\sub S\bigl(C_0(G)\bigr)=\{0\}. 
$$ 
Since this holds for arbitrary $\xi,\eta\in L^2(G)$, 
$S\bigl(B_0(L^2(G))\bigr)=\{0\}$, as claimed. \qed\enddemo

We wish to point out that Proposition 7.8 does not rely on the Hilbert space 
setting: the proof can be modified to yield an analogous result for $L^p(G)$, 
$1<p<\iy$, with the obvious definition of $B_0(L^p(G))$. \ee 
\endexample

\proclaim{Corollary 7.9} Suppose that $E$ is a $W^*$-algebra and 
$\pi(g)\in\Aut(E)$ for every $g\in G$. If $\,\Cal H_{\mu\pi}\ne \{0\}$ and 
$E$ is injective then $\Cal H_{\mu\pi}$ is also injective. 
\endproclaim 
\demo{Proof} We will show that $L_{\mu\pi}$ is injective. Since 
$L^\iy_i(\mu,E)=L^\iy_i(\mu)\oti E$ is injective [40, p.\,120], it suffices 
to construct a norm 1 projection $\Lam$ of $L^\iy_i(\mu,E)$ onto 
$L_{\mu\pi}$. 

Let $\eta$ be any absolutely continuous probability measure on 
$G$. The formula
$$
Jf=\int_G\pi(g)^{-1}f(g)\,\,\eta(dg)
$$
defines a contraction $J$ of $L^\iy(G,E)$ into $E$ such that 
$J[\pi(\cdot)x]_\lam=x$ for every $x\in E$. Define $\hat K=KJ$ where $K$ is 
the projection described in Theorem 7.4. Then $\hat K$ is a contraction of 
$L^\iy(G,E)$ onto $\Cal H_{\mu\pi}$ such that 
$\hat K[\pi(\cdot)x]_\lam =x$ for every $x\in \Cal H_{\mu\pi}$. 

Next, Eq.\,(5.3) defines a contraction, which we will denote again by $R_E$, 
of $L^\iy_i(\mu,E)$ into $L^\iy(G,E)$ $\bigl($in fact, $R_E$ is an isometry 
onto $\Cal H_\mu(E)$, cf. Remark 5.5$\bigr)$. Note that $R_E\vp 
=[\pi(\cdot)R_\pi \vp]_\lam$ for every $\vp\in L_{\mu\pi}$. Put 
$\Lam=R_\pi^{-1}\hat K R_E$. 
\qed\enddemo

Our final result involves the concept of an amenable action of 
a locally compact group on a $W^*$-algebra. The concept 
of an amenable action was first introduced by Zimmer [44] in the context 
of a measure class preserving action of a locally compact second countable 
group on a standard Borel space. 
It was subsequently shown [1] that Zimmer's definition is equivalent to the 
following: Let $\G$ be a locally compact second countable group and 
$\Cal X$ a standard Borel $\G$-space such that the mapping $\G\ti \Cal X\ni 
(g,x)\to gx$ is Borel. Let $\al$ be a $\si$-finite quasiinvariant measure on 
$\Cal X$. Then $\G$ acts on $L^\iy(\Cal X,\al)$ and the resulting 
representation of $\G$ is the adjoint of the natural strongly continuous 
representation of $\G$ in $L^1(\Cal X,\al)$. Consider the tensor 
product $L^\iy(\G)\oti L^\iy(\Cal X,\al)$ equipped with the $\G$-action which 
is the tensor product of the action of $\G$ on $L^\iy(\G)$ by left 
translations and the action on $L^\iy(\Cal X,\al)$. The action of $\G$ on 
$(\Cal X,\al)$ is {\it amenable\/} if there exists a norm 1 equivariant 
projection of $L^\iy(\G)\oti L^\iy(\Cal X,\al)$ onto $L^\iy(\Cal X,\al)$. 
Motivated by the theory of $W^*$-crossed products, this definition was 
extended to the case when $\G$ acts on an arbitrary $W^*$-algebra [2]: Let 
$E$ be a $W^*$-algebra and $\g\:\G\to\Aut(E)$ a representation of 
$\G$ which is the adjoint of a strongly continuous representation of $\G$ 
in $E_*$. The resulting action of $\G$ on $E$ is called {\it amenable} if 
there exists a norm 1 equivariant projection $P$ of $L^\iy(\G)\oti E$ onto 
$E$ where the action of $\G$ on $L^\iy(\G)\oti E$ is the tensor product of 
the action of $\G$ on $L^\iy(\G)$ by left translations and the action $\g$ 
on $E$. We note that by [40, p.\,116] $P$ is necessarily a conditional 
expectation, 
i.e., $P$ preserves positivity and for all $z\in L^\iy(\G)\oti E$ and $x\in 
E$, we have $P\bigl(z(1\oti x)\bigr)=(Pz)x$ and $P\bigl((1\oti 
x)z\bigr)=xPz$. When $\G$ is amenable, then $\g$ is automatically amenable, 
but many actions of nonamenable groups are also amenable. 

Recall that in the case of the classical Choquet-Deny equation in $L^\iy(G)$ 
the action of $G$ on $G$ by left translations gives rise to the 
natural action of $G$ on any $\mu$-boundary. In [44] Zimmer showed that when 
$\mu$ is spread out then this boundary action is always amenable (granted 
that the $\mu$-boundary is a standard Borel $G$-space). In [7] Connes and 
Woods sketched a proof of amenability of the action of $G$ on the space-time 
boundary of a random walk with time-dependent transition probabilities (the 
$\mu$-boundaries discussed here can be viewed as a special case of such 
space-time boundaries). A simple proof of amenability of the action of $G$ on 
boundaries of arbitrary random walks on amenable $G$-spaces given in [14] is 
based on a version of Theorem 7.4. Here we adapt this proof to the setting 
of the Choquet-Deny equation in a $W^*$-algebra $E$. 

Let $\G$ and $\g$ be as described above in the definition of amenablility of 
the action on the $W^*$-algebra $E$ and let us assume that $\g$ and $\pi$ 
commute. Then $\Cal H_{\mu\pi}$ is a $\G$-invariant subspace of $E$ and by 
Corollary 6.10, $\g(g)\r\,\Cal H_{\mu\pi}\in \Aut(\Cal H_{\mu\pi})$ for 
every $g\in \G$. Thus $\G$ acts on $\Cal H_{\mu\pi}$. 

\proclaim{Corollary 7.10}If $\g$ and $\pi$ commute, $\G$ acts amenably 
on the $W^*$-algebra $E$, and if $\,\Cal H_{\mu\pi}\ne\{0\}$, then $\G$ acts 
amenably on $\Cal H_{\mu\pi}$. 
\endproclaim 
\demo{Proof}Let $P\: L^\iy(\G)\oti E\to E$ be the projection appearing in the 
definition of amenability of the action of $\G$ on $E$ and let 
$K$ be the projection described in Theorem 7.4. Then $KP\r 
\bigl(L^\iy(\G)\oti\Cal H_{\mu\pi}\bigr)$ is a norm 1 equivariant 
projection of $L^\iy(\G)\oti \Cal H_{\mu\pi}$ onto $\Cal H_{\mu\pi}$. 
\qed\enddemo

Amenability of the natural action of $\G=G$ on a $\mu$-boundary $(\Cal 
X,\al)$ follows as a special case of Corollary 7.10 because $L^\iy(\Cal 
X,\al)\cong \Cal H_\mu$ and $G$ always acts amenably on $L^\iy(G)$. 

We note that in the case of the ``noncommutative'' Choquet-Deny 
equation considered in Example 2.7, the action of $\G=G$ on $B(L^2(G))$ 
associated with the left regular representation is amenable if and only if 
$G$ is amenable, by [3, Corollaire 3.7]. The trivial example $\mu=\del_e$ 
indicates that, in general, when $G$ fails to be amenable, the action of $G$ 
on the $W^*$-algebra $\Cal H_{\mu\pi}$ of the $\mu$-harmonic operators need 
not be amenable. 

\phantom{oooooooooooooooooo}

\flushpar
{\smc Acknowledgement} 

\flushpar
We are indebted to Dr Christophe Cuny for valuable discussions. 

 \phantom{oooooooooooooooooo}

\flushpar
{\smc Note added in proof}

\flushpar 
Some time after  submission of the present paper 
the authors discovered the article by M\. Izumi  
({\it Non-commutative Poisson boundaries\/},   
Contemp. Math\.  {\bf  347} (2004), 69--81) where 
Example 2.7 is studied,  and the structure formula  
$L^\iy_i(\mu)\ti_{\pi_l}G$  is obtained for the special 
case of a countable discrete group $G$. Izumi asks 
(Problem 4.3) if the latter holds for an arbitrary locally compact 
second countable group.  Our Proposition 6.3 
answers this question in the afirmative .

\enddocument